\documentclass{commat}

\usepackage{enumitem}

\allowdisplaybreaks

\title{%
    On bi-variate poly-Bernoulli polynomials
    }

\author{%
    Claudio Pita-Ruiz
    }

\affiliation{
    \address{Claudio Pita-Ruiz --
    Universidad Panamericana. Facultad de Ingenier\'{\i}a.
Augusto Rodin 498, 
Ciudad de M\'{e}xico, 03920, M\'{e}xico.
        }
    \email{%
    cpita@up.edu.mx
    }
    }

\abstract{%
    We introduce poly-Bernoulli polynomials in two variables by using a
generalization of Stirling numbers of the second kind that we studied in a
previous work. We prove the bi-variate poly-Bernoulli polynomial version of
some known results on standard Bernoulli polynomials, as the addition
formula and the binomial formula. We also prove a result that allows us to
obtain poly-Bernoulli polynomial identities from polynomial identities, and
we use this result to obtain several identities involving products of
poly-Bernoulli and/or standard Bernoulli polynomials. We prove two
generalized recurrences for bi-variate poly-Bernoulli polynomials, and
obtain some corollaries from them.
    }



\VOLUME{31}
\YEAR{2023}
\NUMBER{1}
\firstpage{179}
\DOI{https://doi.org/10.46298/cm.10327}

\begin{paper}

\section{\label{Sec1}Introduction}

Bernoulli numbers are one of the most important mathematical objects that
have been studied by mathemathicians since they appeared in the 18-th
century (see \cite{Ma}). A recent important generalization of the Bernoulli
numbers $B_{n}$ and Bernoulli polynomials $B_{n}(x)$ is about the so-called
poly-Bernoulli numbers and poly-Bernoulli polynomials. Poly-Bernoulli
numbers $B_{n}^{\left( k\right) }$, where $k$ is a given positive integer,
were introduced by M. Kaneko \cite{Ka} in 1997, by means of the generating
function
\begin{equation*}
\frac{\text{\textrm{Li}}_{k}(1-e^{-t})}{1-e^{-t}}=\sum_{n=0}^{\infty
}B_{n}^{\left( k\right) }\frac{t^{n}}{n!},
\end{equation*}
where \textrm{Li}$_{k}(z)=\sum_{j=1}^{\infty }z^{j}/j^{k}$ is the
polylogarithm function. The case $k=1$ corresponds to the standard Bernoulli
numbers $B_{n}$ (except the sign of $B_{1}^{\left( 1\right) }$).
Poly-Bernoulli polynomials $B_{n}^{\left( k\right) }(x)$ can be defined by
the generating function
\begin{equation*}
\frac{\text{\textrm{Li}}_{k}(1-e^{-t})}{1-e^{-t}}e^{-tx}=\sum_{n=0}^{\infty
}B_{n}^{\left( k\right) }(x)\frac{t^{n}}{n!},
\end{equation*}
(see \cite{Ce-Ko}). The case $k=1$ corresponds to $(-1)^{n}B_{n}(x)$, and
the case $x=0$ corresponds to the poly-Bernoulli numbers $B_{n}^{\left(
k\right) }$ mentioned before. Some slightly different definitions of
poly-Bernoulli polynomials $B_{n}^{\left( k\right) }(x)$, with $x$ replaced
by $-x$, and/or with an additional factor $(-1)^{n}$, can be found in some
related papers (see \cite{Ba-Ha}, \cite{Ce-Ko}, \cite{Co-Ca}, \cite{Ha}). In this work we use the
following explicit formula 
\begin{equation}
B_{n}^{\left( k\right) }(x)=\sum_{l=0}^{n}\frac{1}{(l+1)^{k}}
\sum_{j=0}^{l}(-1)^{j}\binom{l}{j}(j+x)^{n},  \label{A}
\end{equation}
as our definition of poly-Bernoulli polynomials (see formula (1.8) in \cite{Ba-Ha}). It is important to mention that the notation
 $B_{n}^{\left(k\right) }(x)$ is also used for a different kind of mathematical objects,
namely, Bernoulli polynomials of $k$-th order (see \cite{Ca1}).

A different generalization of Bernoulli polynomials, studied in the past few
years, is about considering Bernoulli polynomials in several variables $B_{p_{1},\ldots ,p_{t}}(x_{1},\ldots ,x_{t})$, that is, polynomials of
degree $p_{i}$ in the variable $x_{i}$, with 
\[
B_{0, \dotsc, p_{i}, \dotsc, 0}(x_{1},\ldots ,x_{t})
= B_{p_{i}}\left(x_{i}\right) 
\quad \textup{ for each } i \in \{1, \dotsc,t\},
\]
seeking that reasonable generalizations of the known properties in the
one-variable case, remain valid. This kind of work is done in \cite{Shib},
with a flavor of multivariable analysis and working with Jack polynomials. A
different approach is presented in \cite{Shishk} (see also \cite{Bre}, \cite{DiNar}).

In this work we study poly-Bernoulli polynomials in two variables
(bi-variate poly-Bernoulli polynomials). We define the bi-variate
poly-Bernoulli polynomials by using a generalization of Stirling numbers of
the second kind we studied in \cite{Pi}, and then we use the results in \cite{Pi} to obtain results for\ the bi-variate poly-Bernoulli polynomials considered in this work.

We present now the definitions and results in \cite{Pi} that we will use in
the remaining sections.

The generalized Stirling numbers of the second kind (GSN, for short),
denoted as $S_{a_{1},b_{1}}^{\left( a_{2},b_{2},p_{2}\right) }\left(
p_{1},k\right) $, where $a_{j},b_{j}\in \mathbb{C}$, $a_{j}\neq 0$, $j=1,2$,
and $p_{1},p_{2}$ non-negative integers, are defined by means of the
expansion
\begin{equation}
(a_{1}n+b_{1})^{p_{1}}(a_{2}n+b_{2})^{p_{2}}=\sum_{k=0}^{p_{1}+p_{2}}k!S_{a_{1},b_{1}}^{a_{2},b_{2},p_{2}}(p_{1},k)\binom{n}{k},  \label{1.031}
\end{equation}
($S_{a_{1},b_{1}}^{a_{2},b_{2},p_{2}}(p_{1},k)=0$ if $k<0$ or $k>p_{1}+p_{2}$
). An explicit formula for these numbers is 
\begin{equation}
S_{a_{1},b_{1}}^{a_{2},b_{2},p_{2}}(p_{1},k)=\frac{1}{k!} \sum_{j=0}^{k}(-1)^{j}\binom{k}{j} (a_{1}(k-j)+b_{1})^{p_{1}}(a_{2}(k-j)+b_{2})^{p_{2}}.  \label{1.032}
\end{equation}

If $p_{2}=0$, we write the GSN $S_{a,b}^{a_{2},b_{2},0}(p,k)$ as $S_{a,b}(p,k)$. We have 
\begin{equation}
S_{a,b}(p,k)=\frac{1}{k!}\sum_{j=0}^{k}(-1)^{j}\binom{k}{j}\left( a\left(
k-j\right) +b\right) ^{p}.  \label{1.033}
\end{equation}

In the case $a=1,b=0$, the corresponding GSN $S_{1,0}(p,k)=\frac{1}{k!}\sum_{j=0}^{k}(-1)^{j} \binom{k}{j} \left( k-j\right) ^{p}$ are the known Stirling numbers of the second kind. We will refer to them as ``standard
Stirling numbers", and in this case we use the known notation $S(p,k)$.

From (\ref{1.032}) it is clear that $S_{a,b}^{a,b,p_{2}}\left( p_{1},k\right) =S_{a,b}\left( p_{1}+p_{2},k\right) $.
We can see directly from (\ref{1.033}) that 
\begin{eqnarray}
S_{1,1}(p,k) &=&S(p+1,k+1),  \label{1.2} \\
S_{1,2}(p,k) &=&S(p+2,k+2)-S(p+1,k+2).  \label{1.201}
\end{eqnarray}
In this work we will use GSN of the form $S_{1,x_{1}}^{1,x_{2},p_{2}}(p_{1},k)$. Some important facts about the GSN $%
S_{1,x_{1}}^{1,x_{2},p_{2}}(p_{1},k)$ are the following:

\begin{itemize}
\item Some values of the GSN $S_{1,x_{1}}^{1,x_{2},p_{2}}(p_{1},k)$ are
\begin{eqnarray}
S_{1,x_{1}}^{1,x_{2},p_{2}}(p_{1},0) &=&\!\!x_{1}^{p_{1}}\!x_{2}^{p_{2}}\!,
\label{1.3} \\
S_{1,x_{1}}^{1,x_{2},p_{2}}(p_{1},1)
&=&(x_{1}+1)^{p_{1}}\!(x_{2}+1)^{p_{2}}\!-x_{1}^{p_{1}}\!x_{2}^{p_{2}}\!, 
\notag \\
&\vdots&  \notag \\
S_{1,x_{1}}^{1,x_{2},p_{2}}(p_{1},p_{1}+p_{2}) &=&1.  \notag
\end{eqnarray}

\item The GSN $S_{1,x_{1}}^{1,x_{2},p_{2}}(p_{1},k)$ can be written in terms
of the GSN $S_{1,y_{1}}^{1,y_{2},p_{2}}(p_{1},k)$ as follows 
\begin{equation}
S_{1,x_{1}}^{1,x_{2},p_{2}}(p_{1},k)=\sum_{j_{1}=0}^{p_{1}} \sum_{j_{2}=0}^{p_{2}}\binom{p_{1}}{j_{1}} \binom{p_{2}}{j_{2}} (x_{1}-y_{1})^{p_{1}-j_{1}}(x_{2}-y_{2})^{p_{2}-j_{2}}S_{1,y_{1}}^{1,y_{2},j_{2}}(j_{1},k).
\label{1.4}
\end{equation}

\item The GSN $S_{1,x_{1}}^{1,x_{2},p_{2}}(p_{1},k)$ can be written in terms
of standard Stirling numbers as follows
\begin{eqnarray}
k!S_{1,x_{1}}^{1,x_{2},p_{2}}(p_{1},k)
&=&\sum_{j_{1}=0}^{p_{1}}\sum_{j_{2}=0}^{p_{2}}\binom{p_{1}}{j_{1}}\binom{p_{2}}{j_{2}}(x_{1}-m)^{p_{1}-j_{1}}(x_{2}-m)^{p_{2}-j_{2}}  \notag \\
&&\times \sum_{i=0}^{m}\binom{m}{i}(k+i)!S(j_{1}+j_{2},k+i),  \label{1.5}
\end{eqnarray}
where $m$ is an arbitrary non-negative integer.

\item The GSN $S_{1,x_{1}}^{1,x_{2},p_{2}}\left( p_{1},k\right) $ can be
written in terms of standard Stirling numbers as follows
\begin{eqnarray}
S_{1,x_{1}}^{1,x_{2},p_{2}}(p_{1},k)
&=&\sum_{j_{1}=0}^{p_{1}}\sum_{j_{2}=0}^{p_{2}}\binom{p_{1}}{j_{1}}\binom{p_{2}}{j_{2}}(x_{1}-n)^{p_{1}-j_{1}}(x_{2}-n)^{p_{2}-j_{2}}  \notag \\
&&\times \sum_{i=0}^{n-1}(-1)^{i}s(n,n-i)S(j_{1}+j_{2}+n-i,k+n),  \label{1.6}
\end{eqnarray}
where $n$ is an arbitrary positive integer, and $s(\cdot ,\cdot )$ are the
unsigned Stirling numbers of the first kind.

\item The GSN $S_{1,x_{1}}^{1,x_{2},p_{2}}(p_{1},k)$ satisfy the identity
\begin{equation}
S_{1,x_{1}+1}^{1,x_{2}+1,p_{2}}(p_{1},k)=S_{1,x_{1}}^{1,x_{2},p_{2}}(p_{1},k)+(k+1) S_{1,x_{1}}^{1,x_{2},p_{2}}(p_{1},k+1).
\label{1.7}
\end{equation}

\item The GSN $S_{1,x_{1}}^{1,x_{2},p_{2}}(p_{1},k)$ satisfy the recurrence 
\begin{equation}
S_{1,x_{1}}^{1,x_{2},p_{2}}(p_{1},k)=S_{1,x_{1}}^{1,x_{2},p_{2}}(p_{1}-1,k-1) + (k+x_{1})S_{1,x_{1}}^{1,x_{2},p_{2}}(p_{1}-1,k).
\label{1.8}
\end{equation}
\end{itemize}

\section{\label{Sec2}Definitions and preliminary results}

The relation of Bernoulli (numbers and polynomials) with Stirling (numbers
of the second kind) is an old story, that dates back to Worpitsky \cite{W}
(see also \cite[p. 560]{G-K-P} and \cite[p. 5]{Ke}). We have the following
formula for Bernoulli numbers 
\begin{equation}
B_{p}=\sum_{l=0}^{p}S(p,l)\frac{(-1)^{l}l!}{l+1},  \label{2.01}
\end{equation}
and in the case of Bernoulli polynomials we have 
\begin{equation}
B_{p}(x)=\sum_{l=0}^{p}\sum_{j=0}^{p}\binom{p}{j}x^{p-j}S(j,l)\frac{(-1)^{l}l!}{l+1}.  \label{2.02}
\end{equation}

An important observation of formula (\ref{A}) is that poly-Bernoulli
polynomial $B_{p}^{\left( k\right) }(x)$ can be written in terms of the GSN
as
\begin{equation}
B_{p}^{\left( k\right) }(x)=\sum_{l=0}^{p}S_{1,x}(p,l)\frac{(-1)^{l}l!}{(l+1)^{k}}.  \label{2.1}
\end{equation}

The generalization of (\ref{2.1}) to the case of two variables comes through
the GSN: we define poly-Bernoulli polynomial in the variables $x_{1},x_{2}$,
denoted by $B_{p_{1},p_{2}}^{\left( k\right) }(x_{1},x_{2})$, as
\begin{equation}
B_{p_{1},p_{2}}^{\left( k\right)
}(x_{1},x_{2})=\sum_{l=0}^{p_{1}+p_{2}}S_{1,x_{1}}^{1,x_{2},p_{2}}\left(
p_{1},l\right) \frac{(-1)^{l}l!}{(l+1)^{k}}.  \label{2.2}
\end{equation}

If $p_{2}=0$, formula (\ref{2.2}) becomes (\ref{2.1}). By using (\ref{1.5})
we can write $B_{p_{1},p_{2}}^{\left( k\right) }(x_{1},x_{2})$ in terms of
standard Stirling numbers as
\begin{eqnarray}
B_{p_{1},p_{2}}^{\left( k\right) }(x_{1},x_{2})
&=&\sum_{l=0}^{p_{1}+p_{2}}\sum_{j_{1}=0}^{p_{1}}\sum_{j_{2}=0}^{p_{2}} \binom{p_{1}}{j_{1}}\binom{p_{2}}{j_{2}}\left( x_{1}-m\right)^{p_{1}-j_{1}}\left( x_{2}-m\right) ^{p_{2}-j_{2}}  \notag \\
&&\times \sum_{i=0}^{m}\binom{m}{i}S\left( j_{1}+j_{2},l+i\right) \frac{(-1)^{l}(l+i)!}{(l+1)^{k}},  \label{2.5}
\end{eqnarray}
where $m$ is an arbitrary non-negative integer. Similarly, by using (\ref{1.6}) we have that 
\begin{eqnarray}
B_{p_{1},p_{2}}^{\left( k\right) }(x_{1},x_{2})
&=& \sum_{l=0}^{p_{1}+p_{2}}\sum_{j_{1}=0}^{p_{1}}\sum_{j_{2}=0}^{p_{2}} \binom{p_{1}}{j_{1}}\binom{p_{2}}{j_{2}}\left( x_{1}-n\right)^{p_{1}-j_{1}}\left( x_{2}-n\right) ^{p_{2}-j_{2}}  \notag \\
&& \times \sum_{i=0}^{n-1}(-1)^{i}s\left( n,n-i\right) S\left(
j_{1}+j_{2}+n-i,l+n\right) \frac{(-1)^{l}l!}{(l+1)^{k}},  \label{2.6}
\end{eqnarray}
where $n$ is an arbitrary positive integer.

The simplest cases of (\ref{2.5}) and (\ref{2.6}) are 
\begin{equation}
B_{p_{1},p_{2}}^{\left( k\right)}(x_{1},x_{2})  = \sum_{l=0}^{p_{1}+p_{2}}\sum_{j_{1}=0}^{p_{1}} \sum_{j_{2}=0}^{p_{2}}\binom{p_{1}}{j_{1}} \binom{p_{2}}{j_{2}} x_{1}^{p_{1}-j_{1}}x_{2}^{p_{2}-j_{2}}S\left( j_{1}+j_{2},l\right) \frac{(-1)^{l}l!}{(l+1)^{k}},  \label{2.7}
\end{equation}
and
\begin{align}
  B_{p_{1},p_{2}}^{\left( k\right)} & (x_{1},x_{2})   \label{2.8} \\
& = \sum_{l=0}^{p_{1}+p_{2}}\sum_{j_{1}=0}^{p_{1}} \sum_{j_{2}=0}^{p_{2}}\binom{p_{1}}{j_{1}} \binom{p_{2}}{j_{2}} \left(
x_{1}-1\right) ^{p_{1}-j_{1}}\left( x_{2}-1\right) ^{p_{2}-j_{2}} S\left(j_{1}+j_{2}+1,l+1\right) \frac{(-1)^{l}l!}{(l+1)^{k}},  \notag
\end{align}
respectively.

Two examples are the following 
\begin{eqnarray*}
B_{1,1}^{\left( k\right) }(x_{1},x_{2}) &=&\frac{2}{3^{k}}- \frac{1}{2^{k}} \left( x_{1}+x_{2}+1\right) +x_{1}x_{2}, \\
B_{1,2}^{\left( k\right) }(x_{1},x_{2}) &=&\frac{1}{3^{k}}\left(
2x_{1}+4x_{2}+6\right) -\frac{1}{2^{k}}\left( x_{1}\left( 2x_{2}+1\right)
+\left( x_{2}+1\right) ^{2}\right) +x_{1}x_{2}^{2}-\frac{6}{4^{k}}.
\end{eqnarray*}

Clearly we have
\begin{equation}
B_{0,0}^{\left( k\right) }(x_{1},x_{2})=1.  \label{2.9}
\end{equation}

Observe also that
\begin{equation}
B_{p_{1},p_{2}}^{\left( k\right) }\left( x,x\right) =B_{p_{1}+p_{2}}^{\left(
k\right) }(x).  \label{2.10}
\end{equation}

In particular, we have
\begin{equation}
B_{p_{1},p_{2}}^{\left( k\right) }\left( 0,0\right) =B_{p_{1}+p_{2}}^{\left(
k\right) },  \label{2.11}
\end{equation}

From (\ref{1.032}) and (\ref{2.2}) we have
\begin{equation}
B_{p_{1},p_{2}}^{\left( k\right)
}(x_{1},x_{2})=\sum_{l=0}^{p_{1}+p_{2}}\sum_{j=0}^{l}(-1)^{j} \binom{l}{j} \left( l-j+x_{1}\right) ^{p_{1}}\left( l-j+x_{2}\right) ^{p_{2}}\frac{(-1)^{l}}{(l+1)^{k}},  \label{2.13}
\end{equation}
from where we see that 
\begin{eqnarray*}
\frac{\partial }{\partial x_{1}}B_{p_{1},p_{2}}^{\left( k\right)
}(x_{1},x_{2}) &=&p_{1}B_{p_{1}-1,p_{2}}^{\left( k\right) }(x_{1},x_{2}), \\
\frac{\partial }{\partial x_{2}}B_{p_{1},p_{2}}^{\left( k\right)
}(x_{1},x_{2}) &=&p_{2}B_{p_{1},p_{2}-1}^{\left( k\right) }(x_{1},x_{2}).
\end{eqnarray*}

We can use (\ref{1.4}) to write $B_{p_{1},p_{2}}^{\left( k\right)
}(x_{1},x_{2})$ in terms of $B_{j_{1},j_{2}}^{\left( k\right) }\left(
y_{1},y_{2}\right) $, $0\leq j_{i}\leq p_{i}$, $i=1,2$, as 
\begin{equation}
B_{p_{1},p_{2}}^{\left( k\right)}(x_{1},x_{2})= \sum_{j_{1}=0}^{p_{1}} \sum_{j_{2}=0}^{p_{2}} \binom{p_{1}}{j_{1}} \binom{p_{2}}{j_{2}}\left( x_{1}-y_{1}\right) ^{p_{1}-j_{1}} \left(x_{2}-y_{2}\right)^{p_{2}-j_{2}}B_{j_{1},j_{2}}^{\left( k\right) } \left(y_{1},y_{2}\right) ,  \label{2.15}
\end{equation}
that generalizes the known addition formula $B_{p}^{\left( k\right)}(x) = \sum_{j=0}^{p} \binom{p}{j}\left( x-y\right)^{p-j}B_{j}^{\left(k\right) }\left( y\right) $ for one-variable poly-Bernoulli polynomials. In
fact, we have
\begin{align*}
B_{p_{1},p_{2}}^{\left( k\right) }  (x_{1},x_{2})
&=\sum_{l=0}^{p_{1}+p_{2}}S_{1,x_{1}}^{1,x_{2},p_{2}}\left( p_{1},l\right) 
\frac{(-1)^{l}l!}{(l+1)^{k}} \\
&=\sum_{j_{1}=0}^{p_{1}}\sum_{j_{2}=0}^{p_{2}}\binom{p_{1}}{j_{1}}\binom{p_{2}}{j_{2}}(x_{1}-y_{1})^{p_{1}-j_{1}}(x_{2}-y_{2})^{p_{2}-j_{2}} \sum_{l=0}^{j_{1}+j_{2}}S_{1,y_{1}}^{1,y_{2},j_{2}}\left( j_{1},l\right)  \frac{(-1)^{l}l!}{(l+1)^{k}} \\
&=\sum_{j_{1}=0}^{p_{1}}\sum_{j_{2}=0}^{p_{2}}\binom{p_{1}}{j_{1}}\binom{p_{2}}{j_{2}} (x_{1}-y_{1})^{p_{1}-j_{1}}(x_{2}-y_{2})^{p_{2}-j_{2}}B_{j_{1},j_{2}}^{\left( k\right) }(y_{1},y_{2}),
\end{align*}
as claimed. In particular, if we set $y_{1}=y_{2}=y$ in (\ref{2.15}) we
obtain an expression for the bi-variate poly-Bernoulli polynomial $B_{p_{1},p_{2}}^{\left( k\right) }(x_{1},x_{2})$ in terms of one-variable poly-Bernoulli polynomials $B_{j}^{\left( k\right) }\left( y\right) $, $0\leq j\leq p_{1}+p_{2}$, namely
\begin{equation}
B_{p_{1},p_{2}}^{\left( k\right)
}(x_{1},x_{2})=\sum_{j_{1}=0}^{p_{1}} \sum_{j_{2}=0}^{p_{2}} \binom{p_{1}}{j_{1}}\binom{p_{2}}{j_{2}}(x_{1}-y)^{p_{1}-j_{1}}(x_{2}-y)^{p_{2}-j_{2}}B_{j_{1}+j_{2}}^{\left( k\right) }\left( y\right) ,  \label{2.16}
\end{equation}
and then we can write the bi-variate poly-Bernoulli polynomial $B_{p_{1},p_{2}}^{\left( k\right) }(x_{1},x_{2})$ in terms of poly-Bernoulli numbers $B_{j}^{\left( k\right) }$, $0\leq j\leq p_{1}+p_{2}$, as
\begin{equation}
B_{p_{1},p_{2}}^{\left( k\right)}(x_{1},x_{2})=\sum_{j_{1}=0}^{p_{1}}\sum_{j_{2}=0}^{p_{2}}\binom{p_{1}}{j_{1}}\binom{p_{2}}{j_{2}}x_{1}^{p_{1}-j_{1}}x_{2}^{p_{2}-j_{2}}B_{j_{1}+j_{2}}^{\left( k\right) }. \label{2.161}
\end{equation}

The cases $k=0$ and $k=-1$ of (\ref{2.161}) are
\begin{equation}
B_{p_{1},p_{2}}^{\left( 0\right) }(x_{1},x_{2})=\left( x_{1}-1\right)
^{p_{1}}\left( x_{2}-1\right) ^{p_{2}},  \label{2.1611}
\end{equation}
and
\begin{equation}
B_{p_{1},p_{2}}^{(-1)}(x_{1},x_{2})=\left( x_{1}-2\right) ^{p_{1}}\left(
x_{2}-2\right) ^{p_{2}},  \label{2.1612}
\end{equation}
respectively. In fact, according to (\ref{2.2}), formulas (\ref{2.1611}) and
(\ref{2.1612}) are the particular cases $r=0$ and $r=1$ of the identity
\begin{equation}
\sum_{l=0}^{p_{1}+p_{2}}S_{1,x_{1}}^{1,x_{2},p_{2}}\left( p_{1},l\right)
(-1)^{l}\left( l+r\right) !=r!\left( x_{1}-r-1\right) ^{p_{1}}\left(
x_{2}-r-1\right) ^{p_{2}},  \label{2.1615}
\end{equation}
where $r$ is a non-negative integer. We leave the proof of (\ref{2.1615}) to
the reader.

Observe also that (\ref{2.15}) implies that
\begin{equation}
B_{p_{1},p_{2}}^{\left( k\right)
}(x_{1}+1,x_{2}+1)=\sum_{j_{1}=0}^{p_{1}}\sum_{j_{2}=0}^{p_{2}}\binom{p_{1}}{j_{1}}\binom{p_{2}}{j_{2}}B_{j_{1},j_{2}}^{\left( k\right) }(x_{1},x_{2}), \label{2.18}
\end{equation}
which generalizes the known binomial formula for standard Bernoulli
polynomials 
\begin{equation}
B_{p}(x+1)=\sum_{j=0}^{p}\binom{p}{j}B_{j}(x).  \label{2.191}
\end{equation}

If we set $x_{1}=x_{2}=x$ in (\ref{2.15}), we obtain a formula for the
standard poly-Bernoulli polynomial $B_{p_{1}+p_{2}}^{\left( k\right) }(x)$
in terms of the bi-variate poly-Bernoulli polynomials $B_{j_{1},j_{2}}^{\left(
k\right) }(y_{1},y_{2})$, $0\leq j_{i}\leq p_{i}$, $i=1,2$, namely
\begin{equation}
B_{p_{1}+p_{2}}^{\left( k\right)
}(x)=\sum_{j_{1}=0}^{p_{1}}\sum_{j_{2}=0}^{p_{2}}\binom{p_{1}}{j_{1}}\binom{p_{2}}{j_{2}}(x-y_{1})^{p_{1}-j_{1}}(x-y_{2})^{p_{2}-j_{2}}B_{j_{1},j_{2}}^{\left( k\right) }(y_{1},y_{2}).  \label{2.17}
\end{equation}

Some additional simple observations are the following 
\begin{eqnarray}
B_{p_{1},0}^{\left( k\right) }(x_{1},x_{2}) &=&B_{p_{1}}^{\left( k\right)
}(x_{1}),  \label{2.20} \\
B_{0,p_{2}}^{\left( k\right) }(x_{1},x_{2}) &=&B_{p_{2}}^{\left( k\right)
}(x_{2}),  \label{2.201}
\end{eqnarray}
and
\begin{eqnarray}
B_{p_{1},p_{2}}^{\left( k\right) }\left( 0,x_{2}\right)
&=&\sum_{j_{2}=0}^{p_{2}}\binom{p_{2}}{j_{2}} x_{2}^{p_{2}-j_{2}}B_{p_{1}+j_{2}}^{\left( k\right) },  \label{2.21} \\
B_{p_{1},p_{2}}^{\left( k\right) }\left( x_{1},0\right)
&=&\sum_{j_{1}=0}^{p_{1}}\binom{p_{1}}{j_{1}} x_{1}^{p_{1}-j_{1}}B_{j_{1}+p_{2}}^{\left( k\right) }.  \label{2.211}
\end{eqnarray}

\section{\label{Sec3}Some identities}

In this section we obtain some identities involving poly-Bernoulli
polynomials, by using the following result:

\begin{theorem}
\label{Th1}The polynomial identity
\begin{equation}
\sum_{r=0}^{n}a_{n,r}(x+\alpha )^{r}=\sum_{r=0}^{n}b_{n,r}(x+\beta )^{r}. \label{3.1}
\end{equation}
implies the poly-Bernoulli polynomial identity
\begin{equation}
\sum_{r=0}^{n}a_{n,r}B_{r}^{\left( k\right) }(x+\alpha)=\sum_{r=0}^{n}b_{n,r}B_{r}^{\left( k\right) }(x+\beta ).  \label{3.2}
\end{equation}
\end{theorem}

\begin{proof}
Observe that the hypothesis of the polynomial identity (\ref{3.1}), comes
together with the identity of its derivatives:
\begin{equation*}
\sum_{r=0}^{n}\binom{r}{j}a_{n,r}(x+\alpha )^{r-j}=\sum_{r=0}^{n}\binom{r}{j} b_{n,r}(x+\beta )^{r-j}
\end{equation*}
where $j$ is a non-negative integer.

We have
\begin{eqnarray*}
\sum_{r=0}^{n}a_{n,r}B_{r}^{\left( k\right) }(x+\alpha )
&=&\sum_{r=0}^{n}a_{n,r}\sum_{l=0}^{r}S_{1,x+\alpha }(r,l)\frac{(-1)^{l}l!}{(l+1)^{k}} \\
&=&\sum_{r=0}^{n}a_{n,r}\sum_{l=0}^{r}\sum_{j=0}^{r}\binom{r}{j}(x+\alpha
)^{r-j}S(j,l)\frac{(-1)^{l}l!}{(l+1)^{k}} \\
&=&\sum_{l=0}^{n}\sum_{j=0}^{n}\left( \sum_{r=0}^{n} \binom{r}{j}a_{n,r}(x+\alpha )^{r-j}\right) S(j,l) \frac{(-1)^{l}l!}{(l+1)^{k}} \\
&=&\sum_{l=0}^{n}\sum_{j=0}^{n}\left( \sum_{r=0}^{n}\binom{r}{j}b_{n,r}(x+\beta )^{r-j}\right) S(j,l)\frac{(-1)^{l}l!}{(l+1)^{k}} \\
&=&\sum_{r=0}^{n}b_{n,r}\sum_{l=0}^{r}\sum_{j=0}^{r}\binom{r}{j}(x+\beta
)^{r-j}S(j,l)\frac{(-1)^{l}l!}{(l+1)^{k}} \\
&=&\sum_{r=0}^{n}b_{n,r}\sum_{l=0}^{r}S_{1,x+\beta }(r,l)\frac{(-1)^{l}l!}{(l+1)^{k}} \\
&=&\sum_{r=0}^{n}b_{n,r}B_{r}^{\left( k\right) }(x+\beta ),
\end{eqnarray*}
as desired.
\end{proof}

\begin{remark}
The case $k=1$ of Theorem \ref{Th1} is an old result: based on \cite{Ko-Pi},
we obtained Theorem \ref{Th1} in the case $k=1$ and we used it to generate
several identities in \cite{Pi2}.
\end{remark}

For example, by using Theorem \ref{Th1} in the trivial identity $x^{p}=\sum_{j=0}^{p}\binom{p}{j}(x-y)^{j}y^{p-j}$ we obtain the addition
formula for poly-Bernoulli polynomials
\begin{equation}
B_{p}^{\left( k\right) }(x)=\sum_{j=0}^{p}\binom{p}{j}(x-y)^{p-j}B_{j}^{\left( k\right) }(y),  \label{3.3}
\end{equation}
that we can write as
\begin{equation}
\sum_{j=0}^{p} \binom{p}{j}x^{p-j}B_{j}^{\left( k\right) }= \sum_{j=0}^{p} \binom{p}{j} (x-y)^{p-j}B_{j}^{\left( k\right) }(y).  \label{3.4}
\end{equation}

We can use again Theorem \ref{Th1} to obtain from (\ref{3.4}) that
\begin{equation}
\sum_{j=0}^{p}\binom{p}{j}B_{p-j}^{\left( k_{1}\right) }(x)B_{j}^{\left(k\right) }=\sum_{j=0}^{p}\binom{p}{j}B_{p-j}^{\left( k_{1}\right)}(x-y)B_{j}^{\left( k\right) }(y).  \label{3.5}
\end{equation}

Set $y=x$ in (\ref{3.5}) to get the identity
\begin{equation}
\sum_{j=0}^{p}\binom{p}{j}B_{p-j}^{\left( k_{1}\right) }(x)B_{j}^{\left(k\right) } = \sum_{j=0}^{p} 
\binom{p}{j}B_{p-j}^{\left( k_{1}\right)}B_{j}^{\left( k\right) }(x).  \label{3.6}
\end{equation}

If we set $k=1$ in (\ref{3.6}), and replace $x$ by $x+1$, we obtain
\begin{equation}
\sum_{j=0}^{p}\binom{p}{j}B_{p-j}^{\left( k_{1}\right)
}(x+1)B_{j}=\sum_{j=0}^{p}\binom{p}{j}B_{p-j}^{\left( k_{1}\right)
}B_{j}(x+1).  \label{3.7}
\end{equation}

By using that $B_{j}(x+1)=B_{j}(x)+jx^{j-1}$ together with (\ref{3.6}), we
obtain from (\ref{3.7}), after some elementary algebraic steps, the curious
identity 
\begin{equation}
\sum_{j=0}^{p}\binom{p}{j}\left( B_{p-j}^{\left( k_{1}\right)
}(x+1)-B_{p-j}^{\left( k_{1}\right) }(x)\right) B_{j}=pB_{p-1}^{\left(
k_{1}\right) }(x).  \label{3.8}
\end{equation}

From (\ref{2.5}), (\ref{2.6}) and (\ref{2.7}) we have that
\begin{align}
\sum_{j_{1}=0}^{p_{1}}\sum_{j_{2}=0}^{p_{2}}\sum_{l=0}^{j_{1}+j_{2}}  \binom{p_{1}}{j_{1}} & \binom{p_{2}}{j_{2}} x_{1}^{p_{1}-j_{1}}x_{2}^{p_{2}-j_{2}}S(j_{1}+j_{2},l)\frac{(-1)^{l}l!}{(l+1)^{k_{0}}}  \label{3.81} \\
& = \sum_{j_{1}=0}^{p_{1}}\sum_{j_{2}=0}^{p_{2}}\sum_{l=0}^{j_{1}+j_{2}} \binom{p_{1}}{j_{1}}\binom{p_{2}}{j_{2}} (x_{1}-m)^{p_{1}-j_{1}}(x_{2}-m)^{p_{2}-j_{2}} \notag \\
&  \qquad \times  \sum_{i=0}^{m}\binom{m}{i} S(j_{1}+j_{2},l+i)\frac{(-1)^{l}(l+i)!}{(l+1)^{k_{0}}}  \notag \\
& = \sum_{j_{1}=0}^{p_{1}}\sum_{j_{2}=0}^{p_{2}}\sum_{l=0}^{j_{1}+j_{2}} \binom{p_{1}}{j_{1}}\binom{p_{2}}{j_{2}} (x_{1}-n)^{p_{1}-j_{1}}(x_{2}-n)^{p_{2}-j_{2}}  \notag \\
& \qquad \times \sum_{i=0}^{n-1}(-1)^{i}s(n,n-i)S(j_{1}+j_{2}+n-i,l+n)\frac{(-1)^{l}l!}{(l+1)^{k_{0}}}.  \notag
\end{align}
where $m$, $n$ are arbitrary integers, $m\geq 0$, $n>0$. Now we use Theorem \ref{Th1} in (\ref{3.81}) and then set $x_{1}=x_{2}=m+n$, to obtain the
identities
\begin{align}
\sum_{j_{1}=0}^{p_{1}}\sum_{j_{2}=0}^{p_{2}}\sum_{l=0}^{j_{1}+j_{2}} \binom{p_{1}}{j_{1}}\binom{p_{2}}{j_{2}} & B_{p_{1}-j_{1}}^{\left( k_{1}\right)
}(m+n)B_{p_{2}-j_{2}}^{\left( k_{2}\right) }(m+n)S(j_{1}+j_{2},l) \frac{(-1)^{l}l!}{(l+1)^{k_{0}}}  \label{3.82} \\
& = \sum_{j_{1}=0}^{p_{1}}\sum_{j_{2}=0}^{p_{2}}\sum_{l=0}^{j_{1}+j_{2}} \binom{p_{1}}{j_{1}}\binom{p_{2}}{j_{2}}B_{p_{1}-j_{1}}^{\left( k_{1}\right)
}(n)B_{p_{2}-j_{2}}^{\left( k_{2}\right) }(n) \notag  \\ 
& \qquad \times \sum_{i=0}^{m}\binom{m}{i} S(j_{1}+j_{2},l+i)\frac{(-1)^{l}(l+i)!}{(l+1)^{k_{0}}}  \notag \\
& = \sum_{j_{1}=0}^{p_{1}}\sum_{j_{2}=0}^{p_{2}}\sum_{l=0}^{j_{1}+j_{2}} \binom{p_{1}}{j_{1}}\binom{p_{2}}{j_{2}}B_{p_{1}-j_{1}}^{\left( k_{1}\right)
}(m)B_{p_{2}-j_{2}}^{\left( k_{2}\right)}(m) \notag \\ 
& \qquad \times \sum_{i=0}^{n-1}(-1)^{i}s(n,n-i)S(j_{1}+j_{2}+n-i,l+n)\frac{(-1)^{l}l!}{(l+1)^{k_{0}}}.  \notag
\end{align}

From (\ref{2.16}) we see that
\begin{align}
 \sum_{j_{1}=0}^{p_{1}} \sum_{j_{2}=0}^{p_{2}} \binom{p_{1}}{j_{1}} \binom{p_{2}}{j_{2}} & (x_{1}-y)^{p_{1}-j_{1}}(x_{2}-y)^{p_{2}-j_{2}} B_{j_{1}+j_{2}}^{\left( k_{0}\right) }(y)  \label{2.1601} \\
&=\sum_{j_{1}=0}^{p_{1}}\sum_{j_{2}=0}^{p_{2}} \binom{p_{1}}{j_{1}} \binom{p_{2}}{j_{2}}(x_{1}-z)^{p_{1}-j_{1}}(x_{2}-z)^{p_{2}-j_{2}}B_{j_{1}+j_{2}}^{\left( k_{0}\right) }(z).  \notag
\end{align}

We can use Theorem \ref{Th1} to get from (\ref{2.1601}), the identity
\begin{align}
 \sum_{j_{1}=0}^{p_{1}} \sum_{j_{2}=0}^{p_{2}} \binom{p_{1}}{j_{1}} \binom{p_{2}}{j_{2}} & B_{p_{1}-j_{1}}^{\left( k_{1}\right)
}(x_{1}-y) B_{p_{2}-j_{2}}^{\left( k_{2}\right)}(x_{2}-y) B_{j_{1}+j_{2}}^{\left( k_{0}\right) }(y)  \label{2.1602} \\
& = \sum_{j_{1}=0}^{p_{1}} \sum_{j_{2}=0}^{p_{2}} \binom{p_{1}}{j_{1}} \binom{p_{2}}{j_{2}}B_{p_{1}-j_{1}}^{\left( k_{1}\right)} (x_{1}-z)B_{p_{2}-j_{2}}^{\left( k_{2}\right)} (x_{2}-z)B_{j_{1}+j_{2}}^{\left( k_{0}\right) }(z).  \notag
\end{align}

Set $x_{1}=x_{2}=y=x$ to obtain from (\ref{2.1602}) that
\begin{align}
\sum_{j_{1}=0}^{p_{1}}\sum_{j_{2}=0}^{p_{2}} \binom{p_{1}}{j_{1}} \binom{p_{2}}{j_{2}} & B_{p_{1}-j_{1}}^{\left( k_{1}\right) }B_{p_{2}-j_{2}}^{\left(
k_{2}\right) }B_{j_{1}+j_{2}}^{\left( k_{0}\right) }(x)  \label{2.162} \\
& = \sum_{j_{1}=0}^{p_{1}}\sum_{j_{2}=0}^{p_{2}} \binom{p_{1}}{j_{1}} \binom{p_{2}}{j_{2}}B_{p_{1}-j_{1}}^{\left( k_{1}\right)} (x-z)B_{p_{2}-j_{2}}^{\left( k_{2}\right) }(x-z)B_{j_{1}+j_{2}}^{\left( k_{0}\right) }(z).  \notag
\end{align}

With $z=0$, $z=1-\left( q-1\right) x$, and $z=qx-1$, where $q$ is an
arbitrary parameter, we obtain from (\ref{2.162}) the identities
\begin{align}
\sum_{j_{1}=0}^{p_{1}}\sum_{j_{2}=0}^{p_{2}} & \binom{p_{1}}{j_{1}}\binom{p_{2}}{j_{2}} B_{p_{1}-j_{1}}^{\left( k_{1}\right) } B_{p_{2}-j_{2}}^{\left(
k_{2}\right) } B_{j_{1}+j_{2}}^{\left( k_{0}\right) }(x)  \label{2.163} \\
& = \sum_{j_{1}=0}^{p_{1}} \sum_{j_{2}=0}^{p_{2}} \binom{p_{1}}{j_{1}} \binom{p_{2}}{j_{2}}B_{p_{1}-j_{1}}^{\left( k_{1}\right)
}(x)B_{p_{2}-j_{2}}^{\left( k_{2}\right) }(x)B_{j_{1}+j_{2}}^{\left( k_{0}\right) }  \notag \\
& = \sum_{j_{1}=0}^{p_{1}}\sum_{j_{2}=0}^{p_{2}}\binom{p_{1}}{j_{1}} \binom{p_{2}}{j_{2}}B_{p_{1}-j_{1}}^{\left( k_{1}\right)
}(qx-1)B_{p_{2}-j_{2}}^{\left( k_{2}\right) }(qx-1)B_{j_{1}+j_{2}}^{\left( k_{0}\right) }(1-\left( q-1\right) x)  \notag \\
& = \sum_{j_{1}=0}^{p_{1}}\sum_{j_{2}=0}^{p_{2}} \binom{p_{1}}{j_{1}} \binom{p_{2}}{j_{2}}B_{p_{1}-j_{1}}^{\left( k_{1}\right)}(1-(q-1)x)B_{p_{2}-j_{2}}^{\left( k_{2}\right)}(1-(q-1)x)B_{j_{1}+j_{2}}^{\left( k_{0}\right) }(qx-1).  \notag
\end{align}

In the case $q=2$, if some (or all) of the parameters $k_{0},k_{1},k_{2}$
are equal to $1$, we can use the known property $B_{r}\left( 1-x\right)
=(-1)^{r}B_{r}(x)$ to simplify the corresponding expression in (\ref{2.163}). For example, if $k_{0}=k_{1}=k_{2}=1$ and $q=2$, we have the following
identities involving standard Bernoulli numbers and polynomials 

\begin{align}
\sum_{j_{1}=0}^{p_{1}}\sum_{j_{2}=0}^{p_{2}} & \binom{p_{1}}{j_{1}} \binom{p_{2}}{j_{2}} B_{p_{1}-j_{1}}B_{p_{2}-j_{2}} B_{j_{1}+j_{2}}(x)  \label{2.164}
\\
& = \sum_{j_{1}=0}^{p_{1}} \sum_{j_{2}=0}^{p_{2}} \binom{p_{1}}{j_{1}} \binom{p_{2}}{j_{2}}B_{p_{1}-j_{1}}(x)B_{p_{2}-j_{2}}(x)B_{j_{1}+j_{2}}  \notag \\
& = \sum_{j_{1}=0}^{p_{1}} \sum_{j_{2}=0}^{p_{2}} \binom{p_{1}}{j_{1}} \binom{p_{2}}{j_{2}} (-1)^{j_{1}+j_{2}} B_{p_{1}-j_{1}}(2x-1) B_{p_{2}-j_{2}}(2x-1) B_{j_{1}+j_{2}}(x) \notag \\
& = (-1)^{p_{1}+p_{2}} \sum_{j_{1}=0}^{p_{1}} \sum_{j_{2}=0}^{p_{2}} \binom{p_{1}}{j_{1}} \binom{p_{2}}{j_{2}} (-1)^{j_{1}+j_{2}}B_{p_{1}-j_{1}}(x)B_{p_{2}-j_{2}}(x)B_{j_{1}+j_{2}}(2x-1).
\notag
\end{align}

In particular, by setting $x=1$ in (\ref{2.164}), we see that if $p_{1}+p_{2} $ is odd, then
\begin{equation}
\sum_{j_{1}=0}^{p_{1}} \sum_{j_{2}=0}^{p_{2}} \binom{p_{1}}{j_{1}} \binom{p_{2}}{j_{2}}(-1)^{j_{1}+j_{2}}B_{p_{1}-j_{1}}B_{p_{2}-j_{2}}B_{j_{1}+j_{2}}=0.   \label{2.165}
\end{equation}

From (\ref{2.1615}) together with (\ref{1.5}) (with $m=0$) and (\ref{1.6})
(with $n=1$), we have
\begin{align}
\sum_{j_{1}=0}^{p_{1}} & \sum_{j_{2}=0}^{p_{2}} \sum_{l=0}^{j_{1}+j_{2}} \binom{p_{1}}{j_{1}} \binom{p_{2}}{j_{2}} x_{1}^{p_{1}-j_{1}} x_{2}^{p_{2}-j_{2}} S(j_{1}+j_{2},l)(-1)^{l}(l+r)! \label{1.13} \\
& = \sum_{j_{1}=0}^{p_{1}} \sum_{j_{2}=0}^{p_{2}} \sum_{l=0}^{j_{1}+j_{2}} \binom{p_{1}}{j_{1}} \binom{p_{2}}{j_{2}} (x_{1}-1)^{p_{1}-j_{1}}(x_{2}-1)^{p_{2}-j_{2}}S(j_{1}+j_{2}+1,l+1)(-1)^{l}(l+r)! \notag \\
& = r!(x_{1}-r-1)^{p_{1}}(x_{2}-r-1)^{p_{2}},  \notag
\end{align}
and then, applying Theorem \ref{Th1} in (\ref{1.13}) we get 
\begin{align}
& \sum_{j_{1}=0}^{p_{1}} \sum_{j_{2}=0}^{p_{2}} \sum_{l=0}^{j_{1}+j_{2}} \binom{p_{1}}{j_{1}}\binom{p_{2}}{j_{2}} B_{p_{1}-j_{1}}^{\left( k_{1}\right)
}(x_{1})B_{p_{2}-j_{2}}^{\left( k_{2}\right)}(x_{2})S(j_{1}+j_{2},l)(-1)^{l}(l+r)!  \label{1.14} \\
&= \sum_{j_{1}=0}^{p_{1}} \sum_{j_{2}=0}^{p_{2}} \sum_{l=0}^{j_{1}+j_{2}} \binom{p_{1}}{j_{1}} \binom{p_{2}}{j_{2}}B_{p_{1}-j_{1}}^{\left( k_{1}\right)
}(x_{1}-1)B_{p_{2}-j_{2}}^{\left( k_{2}\right)} (x_{2}-1)S(j_{1}+j_{2}+1,l+1)(-1)^{l}(l+r)!  \notag \\
&= r!B_{p_{1}}^{\left( k_{1}\right) }(x_{1}-r-1)B_{p_{2}}^{\left( k_{2}\right) }(x_{2}-r-1),  \notag
\end{align}
where $r$ is an arbitrary non-negative integer.

Now let us consider the difference 
\begin{equation}
B_{p_{1},p_{2}}^{\left( k_{0}\right)} (x_{1}+r,x_{2}+r)-B_{p_{1},p_{2}}^{\left( k_{0}\right) }(x_{1},x_{2}), \label{2.3}
\end{equation}
where $r$ is an arbitrary positive integer. In the case $k_{0}=1$, we know that (\ref{2.3}) is equal to (see \cite{Shib}) 
\begin{equation}
p_{1}\sum_{t=0}^{r-1}(x_{2}+t)^{p_{2}}(x_{1}+t)^{p_{1}-1}+p_{2} \sum_{t=0}^{r-1}(x_{1}+t)^{p_{1}}(x_{2}+t)^{p_{2}-1}.  \label{2.4}
\end{equation}

We can use (\ref{2.16}) to write (\ref{2.3}) as
\begin{align}
\sum_{j_{1}=0}^{p_{1}}\sum_{j_{2}=0}^{p_{2}} & \binom{p_{1}}{j_{1}} \binom{p_{2}}{j_{2}}\left(
(x_{1}+r)^{p_{1}-j_{1}}(x_{2}+r)^{p_{2}-j_{2}}B_{j_{1}+j_{2}}^{\left( k_{0}\right) }-x_{1}^{p_{1}-j_{1}}x_{2}^{p_{2}-j_{2}}B_{j_{1}+j_{2}}^{\left(
k_{0}\right) }\right)  \label{2.401} \\
&=\sum_{j_{1}=0}^{p_{1}}\sum_{j_{2}=0}^{p_{2}} \binom{p_{1}}{j_{1}} \binom{p_{2}}{j_{2}}\left((x_{1}+r-y)^{p_{1}-j_{1}}(x_{2}+r-y)^{p_{2}-j_{2}}B_{j_{1}+j_{2}}^{\left( k_{0}\right) }(y)\right.  \notag \\
& \left. \qquad -(x_{1}-z)^{p_{1}-j_{1}}(x_{2}-z)^{p_{2}-j_{2}}B_{j_{1}+j_{2}}^{\left(k_{0}\right) }(z)\right) ,  \notag
\end{align}
where $y$ and $z$ are arbitrary parameters. If $k_{0}=1$, we have from (\ref{2.4}) and (\ref{2.401}) that
\begin{align}
\sum_{j_{1}=0}^{p_{1}}\sum_{j_{2}=0}^{p_{2}} & \binom{p_{1}}{j_{1}}\binom{p_{2}}{j_{2}}\left( (x_{1}+r)^{p_{1}-j_{1}}(x_{2}+r)^{p_{2}-j_{2}}-x_{1}^{p_{1}-j_{1}}x_{2}^{p_{2}-j_{2}}\right) B_{j_{1}+j_{2}} \label{2.402} \\
& = \sum_{j_{1}=0}^{p_{1}} \sum_{j_{2}=0}^{p_{2}} \binom{p_{1}}{j_{1}} \binom{p_{2}}{j_{2}}\left( (x_{1}+r-y)^{p_{1}-j_{1}}(x_{2}+r-y)^{p_{2}-j_{2}} B_{j_{1}+j_{2}}(y)\right.  \notag \\
& \left. \qquad -(x_{1}-z)^{p_{1}-j_{1}}(x_{2}-z)^{p_{2}-j_{2}} B_{j_{1}+j_{2}}(z)\right)  \notag \\
& =  p_{1}\sum_{t=0}^{r-1}(x_{2}+t)^{p_{2}}(x_{1}+t)^{p_{1}-1}+p_{2} \sum_{t=0}^{r-1}(x_{1}+t)^{p_{1}}(x_{2}+t)^{p_{2}-1}.  \notag
\end{align}

By using Theorem \ref{Th1} in (\ref{2.402}) we get
\begin{align}
\sum_{j_{1}=0}^{p_{1}} \sum_{j_{2}=0}^{p_{2}} & \binom{p_{1}}{j_{1}} \binom{p_{2}}{j_{2}}\left( B_{p_{1}-j_{1}}^{\left( k_{1}\right)} (x_{1}+r) B_{p_{2}-j_{2}}^{\left( k_{2}\right)} (x_{2}+r) - B_{p_{1}-j_{1}}^{\left( k_{1}\right)}(x_{1}) B_{p_{2}-j_{2}}^{\left( k_{2}\right) }(x_{2})\right) B_{j_{1}+j_{2}} \label{2.403} \\
& = \sum_{j_{1}=0}^{p_{1}}\sum_{j_{2}=0}^{p_{2}} \binom{p_{1}}{j_{1}} \binom{p_{2}}{j_{2}}\left( B_{p_{1}-j_{1}}^{\left( k_{1}\right)}(x_{1}+r-y) B_{p_{2}-j_{2}}^{\left( k_{2}\right)} (x_{2}+r-y) B_{j_{1}+j_{2}}(y)\right.  \notag \\
& \left. \qquad - B_{p_{1}-j_{1}}^{\left( k_{1}\right)} (x_{1}-z) B_{p_{2}-j_{2}}^{\left( k_{2}\right)}(x_{2}-z) B_{j_{1}+j_{2}}(z)\right)  \notag \\
& = p_{1} \sum_{t=0}^{r-1} B_{p_{1}-1}^{\left( k_{1}\right)} (x_{1}+t) B_{p_{2}}^{\left( k_{2}\right)} (x_{2}+t)+p_{2}\sum_{t=0}^{r-1} B_{p_{1}}^{\left( k_{1}\right)}(x_{1}+t) B_{p_{2}-1}^{\left( k_{2}\right) }(x_{2}+t).  \notag
\end{align}

Set $\left( y,z\right) =\left( r,r\right) ,\left( r,0\right) ,\left(
0,-r\right) ,\left( 2r,r\right) $ in (\ref{2.403}), to obtain the following
identities involving poly-Bernoulli polynomials, standard Bernoulli
polynomials and Bernoulli numbers

\begin{align}
& \sum_{j_{1}=0}^{p_{1}}\sum_{j_{2}=0}^{p_{2}}  \binom{p_{1}}{j_{1}} \binom{p_{2}}{j_{2}} \left( B_{p_{1}-j_{1}}^{\left( k_{1}\right)} (x_{1}+r) B_{p_{2}-j_{2}}^{\left( k_{2}\right)}(x_{2}+r) - B_{p_{1}-j_{1}}^{\left( k_{1}\right)} (x_{1}) B_{p_{2}-j_{2}}^{\left( k_{2}\right) }(x_{2})\right) B_{j_{1}+j_{2}} \label{1.11} \\
& = \sum_{j_{1}=0}^{p_{1}} \sum_{j_{2}=0}^{p_{2}} \binom{p_{1}}{j_{1}} \binom{p_{2}}{j_{2}} \left( B_{p_{1}-j_{1}}^{\left( k_{1}\right)}(x_{1})B_{p_{2}-j_{2}}^{\left( k_{2}\right)} (x_{2})-B_{p_{1}-j_{1}}^{\left( k_{1}\right)} (x_{1}-r)B_{p_{2}-j_{2}}^{\left( k_{2}\right) }(x_{2}-r)\right) B_{j_{1}+j_{2}}(r)  \notag \\
& = \sum_{j_{1}=0}^{p_{1}}\sum_{j_{2}=0}^{p_{2}}\binom{p_{1}}{j_{1}}\binom{p_{2}}{j_{2}} B_{p_{1}-j_{1}}^{\left( k_{1}\right)}(x_{1})B_{p_{2}-j_{2}}^{\left( k_{2}\right) }(x_{2})\left( B_{j_{1}+j_{2}}(r) - B_{j_{1}+j_{2}}\right)  \notag \\
& = \sum_{j_{1}=0}^{p_{1}}\sum_{j_{2}=0}^{p_{2}} \binom{p_{1}}{j_{1}} \binom{p_{2}}{j_{2}}B_{p_{1}-j_{1}}^{\left( k_{1}\right)} (x_{1}+r) B_{p_{2}-j_{2}}^{\left( k_{2}\right) }(x_{2}+r)\left( B_{j_{1}+j_{2}}-B_{j_{1}+j_{2}}(-r)\right)  \notag \\
& = \sum_{j_{1}=0}^{p_{1}}\sum_{j_{2}=0}^{p_{2}}\binom{p_{1}}{j_{1}}\binom{p_{2}}{j_{2}} B_{p_{1}-j_{1}}^{\left( k_{1}\right)} (x_{1}-r)B_{p_{2}-j_{2}}^{\left( k_{2}\right) }(x_{2}-r)\left( B_{j_{1}+j_{2}}(2r)-B_{j_{1}+j_{2}}(r)\right)  \notag \\
& = p_{1}\sum_{t=0}^{r-1}B_{p_{1}-1}^{\left( k_{1}\right)} (x_{1}+t)B_{p_{2}}^{\left( k_{2}\right)} (x_{2}+t)+p_{2}\sum_{t=0}^{r-1}B_{p_{1}}^{\left( k_{1}\right)
}(x_{1}+t)B_{p_{2}-1}^{\left( k_{2}\right) }(x_{2}+t).  \notag
\end{align}

\section{\label{Sec4}Generalized recurrences}

In this section we show two generalized recurrences for bi-variate
poly-Bernoulli polynomials, and obtain some consequences of them.

\begin{proposition}
\label{Prop4.2}We have
\begin{eqnarray}
\sum_{l=0}^{q}\binom{q}{l}(-x_{1})^{q-l}B_{p_{1}+l,p_{2}}^{\left( k\right)}(x_{1},x_{2}) & = & \sum_{l=0}^{q} \binom{q}{l}(-x_{2})^{q-l}B_{p_{1},p_{2}+l}^{\left( k\right) }(x_{1},x_{2})  \label{4.3} \\
&=&-\sum_{l=0}^{p_{1}+p_{2}}S_{1,x_{1}}^{1,x_{2},p_{2}}(p_{1},l)\frac{(-1)^{l}l!\mathcal{R}_{q-1,k}(l)}{\left( \prod_{i=1}^{q+1}(l+i)\right) ^{k}}.  \notag
\end{eqnarray}
where $\mathcal{R}_{-1,k}(y)=-1$, and for $\mu \geq 0$ the functions $\mathcal{R}_{\mu ,k}(y)$ are defined recursively by 
\begin{equation*}
\mathcal{R}_{\mu ,k}(y)=y(y+\mu +2)^{k}\mathcal{R}_{\mu -1,k}(y)-(y+1)^{k+1} \mathcal{R}_{\mu -1,k}(y+1).
\end{equation*}
\end{proposition}

\begin{proof}
We prove that 
\begin{equation}
\sum_{l=0}^{q}\binom{q}{k}(-x_{1})^{q-l}B_{p_{1}+l,p_{2}}^{\left( k\right)
}(x_{1},x_{2})=-\sum_{l=0}^{p_{1}+p_{2}}S_{1,x_{1}}^{1,x_{2},p_{2}}(p_{1},k) \frac{(-1)^{l}l!\mathcal{R}_{q-1,k}(l)}{\left(
\prod_{i=1}^{q+1}(l+i)\right) ^{k}}  \label{4.4}
\end{equation}
by induction on $q$. The case $q=0$ of (\ref{4.4}) is the definition (\ref{2.2}). Let us suppose formula (\ref{4.4}) is true for a given $q\in \mathbb{N}$. Then 
\begin{align*}
 \sum_{l=0}^{q+1} & \binom{q+1}{l}(-x_{1})^{q+1-l}B_{p_{1}+l,p_{2}}^{\left( k\right) }(x_{1},x_{2}) \\
& = -x_{1}\sum_{l=0}^{q}\binom{q}{l}(-x_{1})^{q-l}B_{p_{1}+l,p_{2}}^{\left( k\right) }(x_{1},x_{2})+\sum_{l=0}^{q}\binom{q}{l}(-x_{1})^{q-l}B_{p_{1}+1+l,p_{2}}^{\left( k\right) }(x_{1},x_{2}) \\
& = x_{1}\sum_{l=0}^{p_{1}}S_{1,x_{1}}^{1,x_{2},p_{2}}(p_{1},l) \frac{(-1)^{l}l!\mathcal{R}_{q-1,k}(l)}{\left( \prod_{i=1}^{q+1}(l+i)\right) ^{k}} - \sum_{l=0}^{p_{1}+1}S_{1,x_{1}}^{1,x_{2},p_{2}}(p_{1}+1,l)\frac{(-1)^{l}l! \mathcal{R}_{q-1,k}(l)}{\left( \prod_{i=1}^{q+1}(l+i)\right) ^{k}}.
\end{align*}
Now we use the recurrence (\ref{1.8}) to write
\begin{align*}
\sum_{l=0}^{q+1}\binom{q+1}{l} & (-x_{1})^{q+1-l}B_{p_{1}+l,p_{2}}^{\left( k\right) }(x_{1},x_{2}) \\
& = x_{1}\sum_{l=0}^{p_{1}}S_{1,x_{1}}^{1,x_{2},p_{2}}(p_{1},l) \frac{(-1)^{l}l! \mathcal{R}_{q-1}(l)}{\left( \prod_{i=1}^{q+1}(l+i) \right) ^{k}} \\
& \qquad -\sum_{l=0}^{p_{1}+1}\left( S_{1,x_{1}}^{1,x_{2},p_{2}}(p_{1},l-1)+(l+x_{1})S_{1,x_{1}}^{1,x_{2},p_{2}}(p_{1},l)\right) 
\frac{(-1)^{l}l!\mathcal{R}_{q-1}(l)}{\left( \prod_{i=1}^{q+1}(l+i)\right)^{k}}.
\end{align*}

Some further simplifications give us 
\begin{align*}
&\sum_{l=0}^{q+1}\binom{q+1}{l}(-x_{1})^{q+1-l}B_{p_{1}+l,p_{2}}^{\left(
k\right) }(x_{1},x_{2}) \\
&=-\sum_{l=0}^{p_{1}+1}\left(
S_{1,x_{1}}^{1,x_{2},p_{2}}(p_{1},l-1)+lS_{1,x_{1}}^{1,x_{2},p_{2}}(p_{1},l) \right) \frac{(-1)^{l}l!\mathcal{R}_{q-1}(l)}{\left(
\prod_{i=1}^{q+1}(l+i)\right) ^{k}} \\
&=\sum_{l=0}^{p_{1}}S_{1,x_{1}}^{1,x_{2},p_{2}}(p_{1},l)\frac{(-1)^{l}(l+1)! \mathcal{R}_{q-1}(l+1)}{\left( \prod_{i=1}^{q+1}(l+i+1)\right) ^{k}} - \sum_{l=0}^{p_{1}}S_{1,x_{1}}^{1,x_{2},p_{2}}(p_{1},l)\frac{(-1)^{l}l!l \mathcal{R}_{q-1}(l)}{\left( \prod_{i=1}^{q+1}(l+i)\right) ^{k}} \\
&=-\sum_{l=0}^{p_{1}}S_{1,x_{1}}^{1,x_{2},p_{2}}\left( p_{1},l\right)
\left( l(l+q+2)^{k}\mathcal{R}_{q-1}(l)-(l+1)^{k+1}\mathcal{R}_{q-1}(l+1)\right) \frac{(-1)^{l}l!}{\left( \prod_{i=1}^{q+2}(l+i)\right)
^{k}} \\
&=-\sum_{l=0}^{p_{1}}S_{1,x_{1}}^{1,x_{2},p_{2}}(p_{1},l)\frac{(-1)^{l}l! \mathcal{R}_{q}(l)}{\left( \prod_{i=1}^{q+2}(l+i)\right) ^{k}},
\end{align*}
as desired. The proof of 
\begin{equation*}
\sum_{l=0}^{q}\binom{q}{l}(-x_{2})^{q-l}B_{p_{1},p_{2}+l}^{\left( k\right)
}(x_{1},x_{2})=-\sum_{l=0}^{p_{1}+p_{2}}S_{1,x_{1}}^{1,x_{2},p_{2}}(p_{1},l) \frac{(-1)^{l}l!\mathcal{R}_{q-1,k}(l)}{\left(
\prod_{i=1}^{q+1}(l+i)\right) ^{k}},
\end{equation*}
is similar.
\end{proof}

For example, we have
\begin{equation*}
\mathcal{R}_{0,k}(y)=(y+1)^{k+1}-y(y+2)^{k},
\end{equation*}
\begin{equation}
\mathcal{R}_{1,k}(y)=(2y+1)(y+1)^{k+1}(y+3)^{k}-y^{2}(y+2)^{k}(y+3)^{k}-(y+1)^{k+1}(y+2)^{k+1},
\label{4.5}
\end{equation}
and then, formula (\ref{4.3}) with $q=1$ is
\begin{align}
-x_{1}B_{p_{1},p_{2}}^{\left( k\right)} (x_{1},x_{2})+ & B_{p_{1}+1,p_{2}}^{\left( k\right) }(x_{1},x_{2})  = -x_{2}B_{p_{1},p_{2}}^{\left( k\right)
}(x_{1},x_{2})+B_{p_{1},p_{2}+1}^{\left( k\right) }(x_{1},x_{2})  \label{4.51} \\
&=-\sum_{l=0}^{p_{1}+p_{2}}S_{1,x_{1}}^{1,x_{2},p_{2}}(p_{1},k)(-1)^{l}l! \left( \frac{l+1}{(l+2)^{k}}-\frac{l}{(l+1)^{k}}\right) ,  \notag
\end{align}
and with $q=2$ is
\begin{align}
x_{1}^{2}B_{p_{1},p_{2}}^{\left( k\right)} & (x_{1},x_{2})-2x_{1}B_{p_{1}+1,p_{2}}^{\left( k\right)} (x_{1},x_{2})+B_{p_{1}+2,p_{2}}^{\left( k\right) }(x_{1},x_{2})  \label{4.6} \\
& = x_{2}^{2}B_{p_{1},p_{2}}^{\left( k\right)}(x_{1},x_{2})-2x_{2}B_{p_{1},p_{2}+1}^{\left( k\right)} (x_{1},x_{2})+B_{p_{1},p_{2}+2}^{\left( k\right) }(x_{1},x_{2})  \notag \\
& = -\sum_{l=0}^{p_{1}+p_{2}}S_{1,x_{1}}^{1,x_{2},p_{2}}(p_{1},l)(-1)^{l}l! \left( \frac{(2l+1)(l+1)}{(l+2)^{k}}- \frac{l^{2}}{(l+1)^{k}}- \frac{(l+1)(l+2)}{(l+3)^{k}}\right) .  \notag
\end{align}

We can write (\ref{4.3}) by using (\ref{2.16}) as 
\begin{align}
\sum_{l=0}^{q} & \binom{q}{l}(-x_{1})^{q-l}\sum_{j_{1}=0}^{p_{1}+l} \sum_{j_{2}=0}^{p_{2}} \binom{p_{1}+l}{j_{1}} \binom{p_{2}}{j_{2}} (x_{1}-y)^{p_{1}+l-j_{1}}(x_{2}-y)^{p_{2}-j_{2}}B_{j_{1}+j_{2}}^{\left( k\right) }(y)  \label{4.61} \\
& = \sum_{l=0}^{q}\binom{q}{l}(-x_{2})^{q-l}\sum_{j_{1}=0}^{p_{1}} \sum_{j_{2}=0}^{p_{2}+l} \binom{p_{1}}{j_{1}} \binom{p_{2}+l}{j_{2}} (x_{1}-z)^{p_{1}-j_{1}}(x_{2}-z)^{p_{2}+l-j_{2}}B_{j_{1}+j_{2}}^{\left( k\right) }(z)  \notag \\
& = -\sum_{l=0}^{p_{1}+p_{2}}S_{1,x_{1}}^{1,x_{2},p_{2}}(p_{1},l) \frac{(-1)^{l}l!\mathcal{R}_{q-1,k}(l)}{\left( \prod_{i=1}^{q+1}(l+i)\right) ^{k}}.  \notag
\end{align}

By setting $y,z=x_{1},x_{2}$ in (\ref{4.61}), we get the identities
\begin{align}
\sum_{l=0}^{q}\binom{q}{l} & (-x_{1})^{q-l} \sum_{j_{2}=0}^{p_{2}}  \binom{p_{2}}{j_{2}}(x_{2}-x_{1})^{p_{2}-j_{2}}B_{p_{1}+l+j_{2}}^{\left( k\right)}(x_{1})  \label{4.62} \\
& = \sum_{l=0}^{q}\binom{q}{l}(-x_{1})^{q-l}\sum_{j_{1}=0}^{p_{1}+l} \binom{p_{1}+l}{j_{1}}(x_{1}-x_{2})^{p_{1}+l-j_{1}}B_{j_{1}+p_{2}}^{\left( k\right)}(x_{2})  \notag \\
& = \sum_{l=0}^{q}\binom{q}{l}(-x_{2})^{q-l}\sum_{j_{2}=0}^{p_{2}+l} \binom{p_{2}+l}{j_{2}}(x_{2}-x_{1})^{p_{2}+l-j_{2}}B_{p_{1}+j_{2}}^{\left( k\right)}(x_{1})  \notag \\
& = \sum_{l=0}^{q}\binom{q}{l}(-x_{2})^{q-l}\sum_{j_{1}=0}^{p_{1}} \binom{p_{1}}{j_{1}}(x_{1}-x_{2})^{p_{1}-j_{1}}B_{j_{1}+p_{2}+l}^{\left( k\right)}(x_{2})  \notag \\
& = -\sum_{j_{1}=0}^{p_{1}}\sum_{j_{2}=0}^{p_{2}}\sum_{l=0}^{j_{1}+j_{2}} \binom{p_{1}}{j_{1}}\binom{p_{2}}{j_{2}} x_{1}^{p_{1}-j_{1}}x_{2}^{p_{2}-j_{2}}S\left( j_{1}+j_{2},l\right) \frac{(-1)^{l}l!\mathcal{R}_{q-1,k}(l)}{\left( \prod_{i=1}^{q+1}(l+i)\right) ^{k}}.  \notag
\end{align}

If we set $x_{2}=0$ in (\ref{4.62}), we get
\begin{align}
\sum_{l=0}^{q}\binom{q}{l}(-x_{1})^{q-l}  \sum_{j_{2}=0}^{p_{2}}  \binom{p_{2}}{j_{2}} & (-x_{1})^{p_{2}-j_{2}}B_{p_{1}+l+j_{2}}^{\left( k\right) }(x_{1}) \label{4.63} \\
& = \sum_{l=0}^{q}\binom{q}{l}(-x_{1})^{q-l}\sum_{j_{1}=0}^{p_{1}+l} \binom{p_{1}+l}{j_{1}}x_{1}^{p_{1}+l-j_{1}}B_{j_{1}+p_{2}}^{\left( k\right) } 
\notag \\
& = \sum_{j_{2}=0}^{p_{2}+q} \binom{p_{2}+q}{j_{2}} (-x_{1})^{p_{2}+q-j_{2}}B_{p_{1}+j_{2}}^{\left( k\right) }(x_{1})  \notag \\
& = \sum_{j_{1}=0}^{p_{1}} \binom{p_{1}}{j_{1}} x_{1}^{p_{1}-j_{1}}B_{j_{1}+p_{2}+q}^{\left( k\right) }  \notag \\
& = -\sum_{l=0}^{p_{1}+p_{2}}S_{1,x_{1}}^{1,0,p_{2}}\left( p_{1},l\right) 
\frac{(-1)^{l}l!\mathcal{R}_{q-1,k}(l)}{\left(
\prod_{i=1}^{q+1}(l+i)\right) ^{k}}.  \notag
\end{align}

The case $q=0,x_{1}=1,k=1$ of (\ref{4.63}) reduces to
\begin{eqnarray}
(-1)^{p_{1}+p_{2}}\sum_{j_{2}=0}^{p_{2}}\binom{p_{2}}{j_{2}}B_{p_{1}+j_{2}}
&=&\sum_{j_{1}=0}^{p_{1}}\binom{p_{1}}{j_{1}}B_{j_{1}+p_{2}}  \label{4.64} \\
&=&\sum_{j_{1}=0}^{p_{1}}\sum_{l=0}^{j_{1}+p_{2}}\binom{p_{1}}{j_{1}} S(j_{1}+p_{2},l)\frac{(-1)^{l}l!}{l+1}.  \notag
\end{eqnarray}

Formula (\ref{4.64}) is the famous Carlitz identity \cite{Ca}. In terms of
bi-variate Bernoulli polynomials, Carlitz identity is written as 
\begin{equation}
B_{p_{1},p_{2}}(1,0)=(-1)^{p_{1}+p_{2}}B_{p_{1},p_{2}}(0,1).  \label{4.65}
\end{equation}

For example, we can use (\ref{4.65}) to write the following version of (\ref{4.3}) in the case $k=1$, when $x_{1}=0,x_{2}=1$
\begin{eqnarray*}
\sum_{l=0}^{q}\binom{q}{l}(-1)^{q-l}B_{p_{1}+l,p_{2}}(1,0)
&=&(-1)^{p_{1}+p_{2}+q}\sum_{l=0}^{q}\binom{q}{l}B_{p_{1}+l,p_{2}}(0,1) \\
&=&B_{p_{1},p_{2}+q}(1,0)=(-1)^{p_{1}+p_{2}+q}B_{p_{1},p_{2}+q}(0,1) \\
&=&-\sum_{l=0}^{p_{1}+p_{2}}S_{1,1}^{1,0,p_{2}}\left( p_{1},l\right) \frac{(-1)^{l}l!\mathcal{R}_{q-1,1}\left( l\right) }{\prod_{i=1}^{q+1}(l+i)},
\end{eqnarray*}
or, explicitly
\begin{align*}
\sum_{l=0}^{q}\binom{q}{l}(-1)^{q-l}\sum_{j_{1}=0}^{p_{1}+l} & \binom{p_{1}+l}{j_{1}}B_{j_{1}+p_{2}}  \\
& =  (-1)^{p_{1}+p_{2}+q}\sum_{l=0}^{q}\binom{q}{l} \sum_{j_{2}=0}^{p_{2}}\binom{p_{2}}{j_{2}}B_{p_{1}+l+j_{2}} \\
&=\sum_{j_{1}=0}^{p_{1}}\binom{p_{1}}{j_{1}} B_{j_{1}+p_{2}+q}=(-1)^{p_{1}+p_{2}+q}\sum_{j_{2}=0}^{p_{2}+q}\binom{p_{2}+q}{j_{2}}B_{p_{1}+j_{2}} \\
&=-\sum_{l=0}^{p_{1}+p_{2}}\sum_{j_{1}=0}^{p_{1}}\binom{p_{1}}{j_{1}} S(j_{1}+p_{2},l)\frac{(-1)^{l}l! \mathcal{R}_{q-1,1}(l)}{\prod_{i=1}^{q+1}(l+i)}.
\end{align*}

It is easy to check that in the case $k=0$, the functions $\mathcal{R}_{\mu
,k}\left( y\right) $ of Proposition \ref{Prop4.2} are $\mathcal{R}_{\mu
,0}(y)=(-1)^{\mu }$. Thus, the case $k=0$ of (\ref{4.3}) is the case $r=0$
of (\ref{2.1615}). Also, in the case $k=-1$, we have that $\mathcal{R}_{\mu
,-1}(y)=\frac{(-1)^{\mu }2^{\mu +1}}{\prod_{i=2}^{\mu +2}(y+i)}$, and then
the case $k=-1$ of (\ref{4.3}) is the case $r=1$ of (\ref{2.1615}).

\begin{proposition}
\label{Prop4.1}For non-negative integers $p_{1},p_{2},q$ we have 
\begin{align}
\sum_{l=0}^{q}  (-1)^{l}B_{p_{1}+l,p_{2}}^{\left( k\right) }(x_{1},x_{2}) \frac{1}{l!} & \frac{d^{l}}{dx_{1}^{l}} \prod_{i=0}^{q-1}(x_{1}+i)  \label{4.1} \\
& =\sum_{l=0}^{q}(-1)^{l}B_{p_{1},p_{2}+l}^{\left( k\right) }(x_{1},x_{2}) \frac{1}{l!}\frac{d^{l}}{dx_{2}^{l}}\prod_{i=0}^{q-1}(x_{2}+i) \notag \\
& = \sum_{l=0}^{p_{1}+p_{2}}S_{1,x_{1}+q}^{1,x_{2}+q,p_{2}}(p_{1},l) \frac{(-1)^{l}(l+q)!}{(l+q+1)^{k}}.  \notag
\end{align}
\end{proposition}

\begin{proof}
We prove that 
\begin{equation}
\sum_{l=0}^{q}(-1)^{l}B_{p_{1}+l,p_{2}}^{\left( k\right) }(x_{1},x_{2})\frac{1}{l!}\frac{d^{l}}{dx_{1}^{l}}\prod_{i=0}^{q-1}(x_{1}+i)=\sum_{l=0}^{p_{1}+p_{2}}S_{1,x_{1}+q}^{1,x_{2}+q,p_{2}}\left( p_{1},l\right) 
\frac{(-1)^{l}\left( l+q\right) !}{\left( l+q+1\right) ^{k}},  \label{4.2}
\end{equation}
by induction on $q$. The case $q=0$ of (\ref{4.2}) is the definition (\ref{2.2}). If we suppose that (\ref{4.2}) is true for a given $q\in \mathbb{N}$, then

\begin{align*}
\sum_{l=0}^{q+1}(-1)^{l} & B_{p_{1}+l,p_{2}}^{\left( k\right) }(x_{1},x_{2}) \frac{1}{l!}\frac{d^{l}}{dx_{1}^{l}}\prod_{i=0}^{q}(x_{1}+i) \\
&=\sum_{l=0}^{q+1}(-1)^{l}B_{p_{1}+l,p_{2}}^{\left( k\right) }(x_{1},x_{2}) \frac{1}{l!}\frac{d^{l}}{dx_{1}^{l}}\left(
(x_{1}+q)\prod_{i=0}^{q-1}(x_{1}+i)\right) \\
&=\sum_{l=0}^{q+1}(-1)^{l}B_{p_{1}+l,p_{2}}^{\left( k\right) }(x_{1},x_{2}) \frac{1}{l!}\left( (x_{1}+q)\frac{d^{l}}{dx_{1}^{l}} \prod_{i=0}^{q-1}(x_{1}+i)+l\frac{d^{l-1}}{dx_{1}^{l-1}} \prod_{i=0}^{q-1}(x_{1}+i)\right) \\
&=(x_{1}+q)\sum_{l=0}^{q}(-1)^{l}B_{p_{1}+l,p_{2}}^{\left( k\right)
}(x_{1},x_{2})\frac{1}{l!}\frac{d^{l}}{dx_{1}^{l}}\prod_{i=0}^{q-1}(x_{1}+i)
\\
& \qquad +\sum_{l=1}^{q+1}(-1)^{l}B_{p_{1}+l,p_{2}}^{\left( k\right) }(x_{1},x_{2}) \frac{1}{\left( l-1\right) !}\frac{d^{l-1}}{dx_{1}^{l-1}} \prod_{i=0}^{q-1}(x_{1}+i) \\
&=(x_{1}+q)\sum_{l=0}^{q}(-1)^{l}B_{p_{1}+l,p_{2}}^{\left( k\right)
}(x_{1},x_{2})\frac{1}{l!}\frac{d^{l}}{dx_{1}^{l}}\prod_{i=0}^{q-1}(x_{1}+i)
\\
& \qquad -\sum_{l=0}^{q}(-1)^{l}B_{p_{1}+1+l,p_{2}}^{\left( k\right) }(x_{1},x_{2}) \frac{1}{l!}\frac{d^{l}}{dx_{1}^{l}}\prod_{i=0}^{q-1}(x_{1}+i) \\
&=(x_{1}+q)\sum_{l=0}^{p_{1}+p_{2}}S_{1,x_{1}+q}^{1,x_{2}+q,p_{2}} \left(p_{1},l\right) \frac{(-1)^{l}\left( l+q\right) !}{\left( l+q+1\right) ^{k}} \\
& \qquad -\sum_{l=0}^{p_{1}+p_{2}+1}S_{1,x_{1}+q}^{1,x_{2}+q,p_{2}} \left(p_{1}+1,l\right) \frac{(-1)^{l}\left( l+q\right) !}{\left( l+q+1\right) ^{k}}.
\end{align*}

Now we use the recurrence (\ref{1.8}) and formula (\ref{1.7}) to write 
\begin{align*}
\sum_{l=0}^{q+1}(-1)^{l} & B_{p_{1}+l,p_{2}}(x_{1},x_{2})\frac{1}{l!} \frac{d^{l}}{dx_{1}^{l}}\prod_{i=0}^{q}(x_{1}+i) \\
&=(x_{1}+q)\sum_{l=0}^{p_{1}+p_{2}}S_{1,x_{1}+q}^{1,x_{2}+q,p_{2}}(p_{1},l) \frac{(-1)^{l}(l+q)!}{(l+q+1)^{k}} \\
& \qquad - \sum_{l=0}^{p_{1}+p_{2}+1} \left( S_{1,x_{1}+q}^{1,x_{2}+q,p_{2}}(p_{1},l-1)+(l+x_{1}+q)S_{1,x_{1}+q}^{1,x_{2}+q,p_{2}}(p_{1},l)\right)  \frac{(-1)^{l}(l+q)!}{(l+q+1)^{k}} \\
&=-\sum_{l=0}^{p_{1}+p_{2}+1} \left( S_{1,x_{1}+q}^{1,x_{2}+q,p_{2}}(p_{1},l-1)+lS_{1,x_{1}+q}^{1,x_{2}+q,p_{2}}(p_{1},l)\right) \frac{(-1)^{l}(l+q)!}{\left( l+q+1\right) ^{k}} \\
& = - \sum_{l=1}^{p_{1}+p_{2}+1}S_{1,x_{1}+q+1}^{1,x_{2}+q+1,p_{2}}(p_{1},l-1) \frac{(-1)^{l}(l+q)!}{\left( l+q+1\right) ^{k}} \\
&= \sum_{l=0}^{p_{1}+p_{2}}S_{1,x_{1}+q+1}^{1,x_{2}+q+1,p_{2}}(p_{1},l) \frac{(-1)^{l}(l+q+1)!}{(l+q+2)^{k}},
\end{align*}
as desired. The proof of 
\begin{equation}
\sum_{l=0}^{q}(-1)^{l}B_{p_{1},p_{2}+l}^{\left( k\right) }(x_{1},x_{2})\frac{1}{l!}\frac{d^{l}}{dx_{2}^{l}}\prod_{i=0}^{q-1}(x_{2}+i) = \sum_{l=0}^{p_{1}+p_{2}}S_{1,x_{1}+q}^{1,x_{2}+q,p_{2}}(p_{1},l)\frac{(-1)^{l}(l+q)!}{(l+q+1)^{k}},  \notag
\end{equation}
is similar.
\end{proof}

Formula (\ref{4.1}) with $x_{1}=0,x_{2}=1$ looks as
\begin{eqnarray}
\sum_{l=0}^{q}(-1)^{l}s(q,l)B_{p_{1}+l,p_{2}}^{\left( k\right) }(0,1)
&=&\sum_{l=0}^{q}(-1)^{l}s(q+1,l+1)B_{p_{1},p_{2}+l}^{\left( k\right) }(0,1)
\label{4.21} \\
&=& \sum_{l=0}^{p_{1}+p_{2}}S_{1,q}^{1,q+1,p_{2}}\left( p_{1},l\right) \frac{(-1)^{l} \left( l+q\right) !}{\left( l+q+1\right) ^{k}},  \notag
\end{eqnarray}
or, explicitly 
\begin{align}
\sum_{l=0}^{q}(-1)^{l} & s(q,l)  \sum_{j_{2}=0}^{p_{2}}\binom{p_{2}}{j_{2}} B_{p_{1}+l+j_{2}}^{\left( k\right) }  
=\sum_{l=0}^{q}(-1)^{l}s(q+1,l+1)\sum_{j_{2}=0}^{p_{2}+l}\binom{p_{2}+l}{j_{2}}B_{p_{1}+j_{2}}^{\left( k\right) } \label{4.22} \\
&=\sum_{l=0}^{p_{1}+p_{2}}\sum_{j_{1}=0}^{p_{1}}\sum_{j_{2}=0}^{p_{2}} \binom{p_{1}}{j_{1}}\binom{p_{2}}{j_{2}} q^{p_{1}-j_{1}}(q+1)^{p_{2}-j_{2}}S(j_{1}+j_{2},l)\frac{(-1)^{l}(l+q)!}{(l+q+1)^{k}}.  \notag
\end{align}

If we set $k=1$ in formula (\ref{4.21}), we can use Carlitz identity (\ref{4.65}) to write the following enriched version of (\ref{4.21}) 
\begin{eqnarray}
\sum_{l=0}^{q}(-1)^{l}s(q,l)B_{p_{1}+l,p_{2}}(0,1)
&=&(-1)^{p_{1}+p_{2}}\sum_{l=0}^{q}s(q,l)B_{p_{1}+l,p_{2}}(1,0)  \label{4.23}
\\
&=&\sum_{l=0}^{q}(-1)^{l}s(q+1,l+1)B_{p_{1},p_{2}+l}(0,1)  \notag \\
&=&(-1)^{p_{1}+p_{2}}\sum_{l=0}^{q}s(q+1,l+1)B_{p_{1},p_{2}+l}(1,0)  \notag
\\
&=&\sum_{l=0}^{p_{1}+p_{2}}S_{1,q}^{1,q+1,p_{2}}(p_{1},l)\frac{(-1)^{l}(l+q)!}{l+q+1},  \notag
\end{eqnarray}
or, explicitly 
\begin{align}
\sum_{l=0}^{q}(-1)^{l}s(q,l) & \sum_{j_{2}=0}^{p_{2}} \binom{p_{2}}{j_{2}} B_{p_{1}+l+j_{2}}  \label{4.24} \\
&=(-1)^{p_{1}+p_{2}}\sum_{l=0}^{q}s(q,l)\sum_{j_{1}=0}^{p_{1}+l} \binom{p_{1}+l}{j_{1}}B_{j_{1}+p_{2}}  \notag \\
&=\sum_{l=0}^{q}(-1)^{l}s(q+1,l+1)\sum_{j_{2}=0}^{p_{2}+l} \binom{p_{2}+l}{j_{2}}B_{p_{1}+j_{2}}  \notag \\
&=(-1)^{p_{1}+p_{2}}\sum_{l=0}^{q}s(q+1,l+1)\sum_{j_{1}=0}^{p_{1}} \binom{p_{1}}{j_{1}}B_{j_{1}+p_{2}+l}  \notag \\
&=\sum_{j_{1}=0}^{p_{1}}\sum_{j_{2}=0}^{p_{2}}\sum_{l=0}^{j_{1}+j_{2}} \binom{p_{1}}{j_{1}} \binom{p_{2}}{j_{2}}q^{p_{1}-j_{1}}(q+1)^{p_{2}-j_{2}}S(j_{1}+j_{2},l) \frac{(-1)^{l}(l+q)!}{l+q+1}.  \notag
\end{align}

To end this section, let us consider the case $q=1$ of the first two lines
of (\ref{4.1}). That is
\begin{equation}
B_{p_{1},p_{2}+1}^{\left( k\right) }(x_{1},x_{2})-B_{p_{1}+1,p_{2}}^{\left(
k\right) }(x_{1},x_{2})=(x_{2}-x_{1})B_{p_{1},p_{2}}^{\left( k\right)
}(x_{1},x_{2}).  \label{4.7}
\end{equation}

Formula (\ref{4.7}) is the first step of two results contained in the
following proposition.

\begin{proposition}
We have the following identities:
\begin{enumerate}[label=\alph*)]

\item
\begin{equation}
\sum_{j=0}^{q}\binom{q}{j}(-1)^{j}B_{p_{1}+j,p_{2}+q-j}^{\left( k\right)
}(x_{1},x_{2})=(x_{2}-x_{1})^{q}B_{p_{1},p_{2}}^{\left( k\right)
}(x_{1},x_{2}).  \label{4.8}
\end{equation}

\item
\begin{equation}
\sum_{j=0}^{q}\binom{q}{j}(x_{2}-x_{1})^{j}B_{p_{1}+q-j,p_{2}}^{\left(
k\right) }(x_{1},x_{2})=B_{p_{1},p_{2}+q}^{\left( k\right) }(x_{1},x_{2}).
\label{4.9}
\end{equation}
\end{enumerate}
\end{proposition}

\begin{proof}
Let us prove (\ref{4.8}) by induction on $q$. The case $q=0$ is a trivial
identity. Let us suppose that (\ref{4.8}) is true for a given $q\in \mathbb{N%
}$. Then 
\begin{align*}
\sum_{j=0}^{q+1} & \binom{q+1}{j}(-1)^{j}B_{p_{1}+j,p_{2}+q+1-j}^{\left(
k\right) }(x_{1},x_{2}) \\
&=\sum_{j=0}^{q+1}\left( \binom{q}{j}+\binom{q}{j-1}\right)
(-1)^{j}B_{p_{1}+j,p_{2}+q+1-j}^{\left( k\right) }(x_{1},x_{2}) \\
&=\sum_{j=0}^{q}\binom{q}{j}(-1)^{j}B_{p_{1}+j,p_{2}+q+1-j}^{\left(
k\right) }(x_{1},x_{2})+\sum_{j=0}^{q}\binom{q}{j}(-1)^{j+1}B_{p_{1}+1+j,p_{2}+q-j}^{\left( k\right) }(x_{1},x_{2}) \\
&=(x_{2}-x_{1})^{q}B_{p_{1},p_{2}+1}^{\left( k\right)
}(x_{1},x_{2})-(x_{2}-x_{1})^{q}B_{p_{1}+1,p_{2}}^{\left( k\right)
}(x_{1},x_{2}) \\
&=(x_{2}-x_{1})^{q}\left( B_{p_{1},p_{2}+1}^{\left( k\right)
}(x_{1},x_{2})-B_{p_{1}+1,p_{2}}^{\left( k\right) }(x_{1},x_{2})\right) \\
&=(x_{2}-x_{1})^{q+1}B_{p_{1},p_{2}}^{\left( k\right) }(x_{1},x_{2}),
\end{align*}
as desired. In the last step we used (\ref{4.7}).

Now let us prove (\ref{4.9}). Again we proceed by induction on $q$. The case 
$q=0$ is a trivial identity. Let us suppose that (\ref{4.9}) is true for a
given $q\in \mathbb{N}$. Then
\begin{align*}
\sum_{j=0}^{q+1} & \binom{q+1}{j}  (x_{2}-x_{1})^{j}B_{p_{1}+q+1-j,p_{2}}^{\left( k\right) }(x_{1},x_{2}) \\
&=\sum_{j=0}^{q+1}\left( \binom{q}{j}+\binom{q}{j-1}\right)
(x_{2}-x_{1})^{j}B_{p_{1}+q+1-j,p_{2}}^{\left( k\right) }(x_{1},x_{2}) \\
&=\sum_{j=0}^{q}\binom{q}{j}(x_{2}-x_{1})^{j}B_{p_{1}+q+1-j,p_{2}}^{\left(
k\right) }(x_{1},x_{2})+\sum_{j=0}^{q}\binom{q}{j} (x_{2}-x_{1})^{j+1}B_{p_{1}+q-j,p_{2}}^{\left( k\right) }(x_{1},x_{2}) \\
&= B_{p_{1} + 1,p_{2}+q}^{\left( k\right)}(x_{1},x_{2})+(x_{2}-x_{1})B_{p_{1},p_{2}+q}^{\left( k\right) }(x_{1},x_{2})
\\
&=B_{p_{1},p_{2}+q+1}^{\left( k\right) }(x_{1},x_{2}),
\end{align*}
as desired. We used (\ref{4.7}) (with $p_{2}$ replaced by $p_{2}+q$) in the
last step.
\end{proof}

The case $p_{1}=0$ of (\ref{4.9}) is 
\begin{equation}
\sum_{j=0}^{q}\binom{q}{j}(x_{2}-x_{1})^{j}B_{q-j,p_{2}}^{\left( k\right)
}(x_{1},x_{2})=B_{p_{2}+q}^{\left( k\right) }(x_{2})\text{.}  \label{4.10}
\end{equation}

The case $p_{2}=0$ of (\ref{4.10}) is the addition formula for standard
poly-Bernoulli polynomials, namely 
\begin{equation*}
\sum_{j=0}^{q}\binom{q}{j}(x_{2}-x_{1})^{q-j}B_{j}^{\left( k\right)
}(x_{1})=B_{q}^{\left( k\right) }(x_{2}).
\end{equation*}

The case $p_{1}=p_{2}=0$ of (\ref{4.8}) is 
\begin{equation}
\sum_{j=0}^{q}\binom{q}{j}(-1)^{j}B_{j,q-j}^{\left( k\right)
}(x_{1},x_{2})=(x_{2}-x_{1})^{q}.  \label{4.11}
\end{equation}

As a final comment, we mention that by considering the GSN $S_{1,x_{1}}^{(1,x_{2},p_{2}),\ldots ,(1,x_{n},p_{n})}\left( p_{1},k\right) $
involved in the expansion
\begin{equation}
(m+x_{1})^{p_{1}}\cdots (m+x_{n})^{p_{n}}=\sum_{l=0}^{p_{1}+\cdots
+p_{n}}l!S_{1,x_{1}}^{(1,x_{2},p_{2}),\ldots ,(1,x_{n},p_{n})}\left(
p_{1},l\right) \binom{m}{l},
\end{equation}
where $p_{1},p_{2},\ldots ,p_{n}$ are non-negative integers given, one can
define poly-Bernoulli polynomials in $n$ variables $x_{1},\ldots ,x_{n}$,
denoted as $B_{p_{1},\ldots ,p_{n}}^{(k)}(x_{1},\ldots ,x_{n})$, as 
\begin{equation*}
B_{p_{1},\ldots ,p_{n}}^{\left( k\right) }(x_{1},\ldots
,x_{n})=\sum_{l=0}^{p_{1}+\cdots +p_{n}}S_{1,x_{1}}^{(1,x_{2},p_{2}),\ldots
,(1,x_{n},p_{n})}\left( p_{1},l\right) \frac{(-1)^{l}l!}{(l+1)^{k}},
\end{equation*}
or explicitly as 
\begin{align*}
B_{p_{1},\ldots ,p_{n}}^{\left( k\right) } & (x_{1},\ldots ,x_{n}) \\
&=\sum_{l=0}^{p_{1}+\cdots +p_{n}}\sum_{j_{1}=0}^{p_{1}}\cdots
\sum_{j_{n}=0}^{p_{n}}\binom{p_{1}}{j_{1}}\cdots \binom{p_{n}}{j_{n}}x_{1}^{p_{1}-j_{1}}\cdots x_{n}^{p_{n}-j_{n}}S\left( j_{1}+\cdots
+j_{n},l\right) \frac{(-1)^{l}l!}{(l+1)^{k}} \\
&=\sum_{j_{1}=0}^{p_{1}}\cdots \sum_{j_{n}=0}^{p_{n}}\binom{p_{1}}{j_{1}}\cdots \binom{p_{n}}{j_{n}}x_{1}^{p_{1}-j_{1}}\cdots
x_{n}^{p_{n}-j_{n}}B_{j_{1}+\cdots +j_{n}}^{\left( k\right) }.
\end{align*}

In this more general setting we have natural generalizations of results (\ref{4.3}), (\ref{4.1}), (\ref{4.8}) and (\ref{4.9}). We show the corresponding
results in the case of poly-Bernoulli polynomials in 3 variables: 
\begin{equation*}
B_{p_{1},p_{2},p_{3}}^{(k)}(x_{1},x_{2},x_{3})=\sum_{j_{1}=0}^{p_{1}}\sum_{j_{2}=0}^{p_{2}}\sum_{j_{3}=0}^{p_{3}}\binom{p_{1}}{j_{1}}\binom{p_{2}}{j_{2}}\binom{p_{3}}{j_{3}}x_{1}^{p_{1}-j_{1}}x_{2}^{p_{2}-j_{2}}x_{3}^{p_{3}-j_{3}}B_{j_{1}+j_{2}+j_{3}}^{\left( k\right) },
\end{equation*}

\begin{enumerate}[label=(\alph*)]
\item (See (\ref{4.3})) We have the generalized recurrences 
\begin{align*}
\sum_{l=0}^{q}\binom{q}{l} & (-x_{1})^{q-l}B_{p_{1}+l,p_{2},p_{3}}^{\left(k\right) }(x_{1},x_{2},x_{3}) \\
&=\sum_{l=0}^{q} \binom{q}{l}(-x_{2})^{q-l}B_{p_{1},p_{2}+l,p_{3}}^{\left(k\right) }(x_{1},x_{2},x_{3}) \\
&=\sum_{l=0}^{q} \binom{q}{l}(-x_{3})^{q-l}B_{p_{1},p_{2},p_{3}+l}^{\left(k\right) }(x_{1},x_{2},x_{3}) \\
&=- \sum_{l=0}^{p_{1}+p_{2}+p_{3}} \sum_{j_{1}=0}^{p_{1}} \sum_{j_{2}=0}^{p_{2}}\sum_{j_{3}=0}^{p_{3}} \left( \binom{p_{1}}{j_{1}}\binom{p_{2}}{j_{2}} \binom{p_{3}}{j_{3}} x_{1}^{p_{1}-j_{1}}x_{2}^{p_{2}-j_{2}}x_{3}^{p_{3}-j_{3}} \right.\\
& \hspace{6cm} \left. \times S(j_{1}+j_{2}+j_{3},l)\frac{(-1)^{l}l!\mathcal{R}_{q-1,k}(l)}{\left(\prod_{i=1}^{q+1}(l+i)\right) ^{k}} \right) ,
\end{align*}
where $\mathcal{R}_{q-1,k}(l)$ is defined in Proposition \ref{Prop4.2}.

\item (See (\ref{4.1})) We have the generalized recurrences
\begin{align*}
\sum_{l=0}^{q} & (-1)^{l}  B_{p_{1}+l,p_{2},p_{3}}^{\left( k\right)}(x_{1},x_{2},x_{3})\frac{1}{l!}  \frac{d^{l}}{dx_{1}^{l}} \prod_{i=0}^{q-1}(x_{1}+i)  \\
&=\sum_{l=0}^{q}(-1)^{l}B_{p_{1},p_{2}+l,p_{3}}^{\left( k\right)}(x_{1},x_{2},x_{3})\frac{1}{l!}  \frac{d^{l}}{dx_{2}^{l}} \prod_{i=0}^{q-1}(x_{2}+i)  \\
&=\sum_{l=0}^{q}(-1)^{l}B_{p_{1},p_{2},p_{3}+l}^{\left( k\right)}(x_{1},x_{2},x_{3})\frac{1}{l!}  \frac{d^{l}}{dx_{3}^{l}} \prod_{i=0}^{q-1}(x_{3}+i)  \\
&=\sum_{j_{1}=0}^{p_{1}}\sum_{j_{2}=0}^{p_{2}}\sum_{j_{3}=0}^{p_{3}} \left( \binom{p_{1}}{j_{1}}\binom{p_{2}}{j_{2}}\binom{p_{3}}{j_{3}}   (x_{1}+q)^{p_{1}-j_{1}}(x_{2}+q)^{p_{2}-j_{2}}  (x_{3}+q)^{p_{3}-j_{3}}  \right. \\
&   \hspace{7cm} \left. \times S(j_{1}+j_{2}+j_{3},l) \frac{(-1)^{l}(l+q)!}{(l+q+1)^{k}} \right).
\end{align*}

\item (See (\ref{4.8})) We have the identities
\begin{eqnarray*}
\sum_{j=0}^{q}\binom{q}{j}\left( -1\right)
^{j}B_{p_{1},p_{2}+j,p_{3}+q-j}^{\left( k\right) }\left(
x_{1},x_{2},x_{3}\right) &=&\left( x_{3}-x_{2}\right)
^{q}B_{p_{1},p_{2},p_{3}}^{\left( k\right) }\left( x_{1},x_{2},x_{3}\right) ,
\\
\sum_{j=0}^{q}\binom{q}{j}\left( -1\right)
^{j}B_{p_{1}+j,p_{2},p_{3}+q-j}^{\left( k\right) }\left(
x_{1},x_{2},x_{3}\right) &=&\left( x_{3}-x_{1}\right)
^{q}B_{p_{1},p_{2},p_{3}}^{\left( k\right) }\left( x_{1},x_{2},x_{3}\right) ,
\\
\sum_{j=0}^{q}\binom{q}{j}\left( -1\right)
^{j}B_{p_{1}+j,p_{2}+q-j,p_{3}}^{\left( k\right) }\left(
x_{1},x_{2},x_{3}\right) &=&\left( x_{2}-x_{1}\right)
^{q}B_{p_{1},p_{2},p_{3}}^{\left( k\right) }\left( x_{1},x_{2},x_{3}\right) .
\end{eqnarray*}

\item (See (\ref{4.9})) We have the identities
\begin{align*}
\sum_{j=0}^{q}\binom{q}{j} & \left( x_{1}-x_{2}\right)^{j}B_{p_{1},p_{2}+q-j,p_{3}}^{\left( k\right) } \left(
x_{1},x_{2},x_{3}\right) \\
& =  \sum_{j=0}^{q} \binom{q}{j}\left( x_{1}-x_{3}\right)^{j}B_{p_{1},p_{2},p_{3}+q-j}^{\left( k\right) }\left(
x_{1},x_{2},x_{3}\right) = B_{p_{1}+q,p_{2},p_{3}}^{\left( k\right) } \left( x_{1},x_{2},x_{3}\right).
\end{align*}

\begin{align*}
\sum_{j=0}^{q}\binom{q}{j} & \left( x_{2}-x_{1}\right)
^{j}B_{p_{1}+q-j,p_{2},p_{3}}^{\left( k\right) }\left(
x_{1},x_{2},x_{3}\right) \\
&= \sum_{j=0}^{q} \binom{q}{j} \left( x_{2}-x_{3}\right)^{j}B_{p_{1},p_{2},p_{3}+q-j}^{\left( k\right) } \left(x_{1},x_{2},x_{3}\right) = B_{p_{1},p_{2}+q,p_{3}}^{\left( k\right) }\left( x_{1},x_{2},x_{3}\right).
\end{align*}

\begin{align*}
\sum_{j=0}^{q}\binom{q}{j} & \left( x_{3}-x_{1}\right)^{j}B_{p_{1}+q-j,p_{2},p_{3}}^{\left( k\right) }\left(
x_{1},x_{2},x_{3}\right) \\
& = \sum_{j=0}^{q} \binom{q}{j}\left( x_{3}-x_{2}\right)^{j} B_{p_{1},p_{2}+q-j,p_{3}}^{\left( k\right) }\left(
x_{1},x_{2},x_{3}\right) = B_{p_{1},p_{2},p_{3}+q}^{\left( k\right) } \left( x_{1},x_{2},x_{3}\right).
\end{align*}

\end{enumerate}

\subsection*{Acknowledgments}

I thank the anonymous referee for his/her helpful observations, that
certainly contributed to improve the final version of this work.


\EditInfo{October 08, 2020}{November 28, 2021}{Karl Dilcher}

\end{paper}